\documentclass[onefignum,onetabnum]{siamart220329}
\usepackage{import}
\usepackage{ownPackage}

\usepackage{todonotes}

\headers{Deflation for symmetric saddle point systems}{A. Dumitrasc, C. Kruse and U. R\" ude}

\begin{document}

	\title{Deflation for the off-diagonal block in symmetric saddle point systems
		\thanks{\funding{This research was supported by the Bavarian Academic Center for Central, Eastern and Southeastern Europe (BAYHOST).}}}
	
	\author{A. Dumitrasc\thanks{Chair for Computer Science 10 - System Simulation, Friedrich-Alexander University Erlangen-Nuremberg, 91058, Erlangen, Germany 
			(\email{andrei.dumitrasc@fau.de}, \email{ulrich.ruede@fau.de}).}		
		\and C. Kruse\thanks{CERFACS, 31057, Toulouse, France
			(\email{carola.kruse@cerfacs.fr}).}	
		\and U. R\" ude\footnotemark[1] \footnotemark[2]
	}
	
	\maketitle

	\begin{abstract}
		Deflation techniques are typically used to shift isolated clusters of small eigenvalues in order to obtain a tighter distribution and a smaller condition number. Such changes induce a positive effect in the convergence behavior of Krylov subspace methods, which are among the most popular iterative solvers for large sparse linear systems. We develop a deflation strategy for symmetric saddle point matrices by taking advantage of their underlying block structure. The vectors used for deflation come from an elliptic singular value decomposition relying on the generalized Golub-Kahan bidiagonalization process. The block targeted by deflation is the off-diagonal one since it features a problematic singular value distribution for certain applications. One example is the Stokes flow in elongated channels, where the off-diagonal block has several small, isolated singular values, depending on the length of the channel. Applying deflation to specific parts of the saddle point system is important when using solvers such as CRAIG, which operates on individual blocks rather than the whole system. The theory is developed by extending the existing framework for deflating square matrices before applying a Krylov subspace method like MINRES. Numerical experiments confirm the merits of our strategy and lead to interesting questions about using approximate vectors for deflation.
	\end{abstract}

	

	 \section{Introduction}
	 \label{sec:intro_CK}

	 Advanced simulation tasks often require the numerical solution
	 of partial differential equations (PDE) with fine resolution 
	 such that the discretization leads to very large, 
	 sparse linear systems. For many practical applications,
	 the computational cost of the linear solvers is among the critical factors,
	 since it can be dominating and a bottleneck in the overall cost of the simulation and thus its speed. While direct methods with their generality and high accuracy often represent the first choice, they may not be efficient or their application may even become unfeasible due to their high computational complexity. As alternative, iterative solvers that successively improve an initial guess of the solution vector can be used. These solvers can be significantly cheaper, also as the iteration can be stopped at any time, 
	 when only a moderately precise solution is needed. However, they rarely work well in a black-box manner and have to be tailored to the given matrix properties that may be hard to identify.
	 
	 One widely used class of iterative solvers are 
	 Krylov subspace methods.
	 In the case of algorithms designed for symmetric problems, their speed of convergence is linked to the conditioning or, more precisely, the spectral properties of the system matrix 
	 \cite{saad2003iterative, barrett1994templates}. 
	 To obtain fast convergence, it is desirable to have clustered eigenvalues that are bounded away from zero. 
	 When this is not the case, the method converges slowly or can even stagnate \cite{liesen2004convergence}. The general strategy is then to help the iterative solver by manipulating the linear system, for example by preconditioning.
	 Preconditioning aims to improve the spectral properties of the matrix, such that the linear system becomes easier to  
	 solve \cite{barrett1994templates}. 
	 Depending on how much one wishes to invest, the preconditioner can be as simple as a diagonal matrix, or as sophisticated as a matrix representing the eigenvectors associated with the smallest eigenvalues. An alternative technique to preconditioning is known as deflation.
	 It is similar, in the sense that it affects the spectrum of the matrix and makes the system easier to solve, by leading to a better clustering of the eigenvalues. Unlike preconditioning, it has a more targeted effect, shifting only the outliers such that they are either zero or closer to the main cluster. The rest of the spectrum is left unchanged. Note that this is not the only definition associated with the idea of deflation. For a concise but informative chronological overview of the developments related to deflation we refer to \cite{gutknecht2012spectral}. A notable work describing deflation and augmentation methods under the wider umbrella of recycling is the more recent survey \cite{soodhalter2020survey}, where significant attention is also given to sequences of (shifted) systems. An example of such a Krylov method with recycling is RMINRES, used successfully in an application of electrical impedance tomography (see \cite{mello2010recycling}).
	 
	 In this paper, we focus on iterative solvers for matrices in saddle point form, as they appear, for example, in the discretization of PDEs in incompressible fluid flow, which is a longstanding research focus in our group \cite{KrDaTaArRu2020,KrSoArTaRu2020,KrSoArTaRu2021,darrigrand2022inexact}. In particular, we focus on fluid flow in a thin-elongated channel domain, where the ratio of length and the height is large, as described by the Stokes equation
	 \begin{align*}
	 	- \Delta \Vec{u}  + \nabla p &= \Vec{f} \\
	 	\nabla \cdot  \Vec{u} &=0.
	 \end{align*}
	 The exact definition, boundary conditions and special properties of our test problem are given in \Cref{sec:testPb}. Here, we only want to sketch the behavior and difficulty of the chosen iterative solver and the aim of our work. When discretizing the above equation by, e.g., finite differences or a stable class of finite elements, we obtain a linear system with a matrix of $2 \times 2$ block structure 
	 
	 \begin{equation}\label{eq:sad_po}
	 	\left[
	 	\begin{matrix}
	 		\Wm & \Am \\
	 		\Am ^{T} & \mZ 
	 	\end{matrix}
	 	\right]
	 	\left[
	 	\begin{matrix}
	 		\uv \\
	 		\pv 
	 	\end{matrix}
	 	\right]=
	 	\left[
	 	\begin{matrix}
	 		\gvv \\
	 		\rv 
	 	\end{matrix}
	 	\right], 
	 \end{equation}
	 with $\Wm \in \bR^{m\times m}$ being symmetric and positive definite (SPD) and $\Am \in \bR^{m \times n}, (n<m)$ having full column rank. The vectors are $\uv \in \bR^{m}, \pv \in \bR^{n}, \gvv \in \bR^{m}$ and $ \rv  \in \bR^{n}$, with $\uv $ and $ \pv$ representing velocity and pressure, respectively. This saddle point system matrix is indefinite and thus has positive and negative eigenvalues. An overview about solution methods, including \textit{segregated methods} that exploit the block structure of the system and \textit{monolithic methods} that are applied to the system as a whole, is given in \cite{benzi2005numerical}. While it is known that, e.g., monolithic multigrid methods can achieve excellent performance for special saddle point problems \cite{kohl2020textbook}, a block version offers 
	 the possibility to re-use existing solvers. In the following work, we are interested in solving the above linear system by Krylov subspace methods. Since the system is symmetric, we could apply as a first choice the well-known \textit{minimum residual} (MINRES) method monolithically. Its convergence curve is characterized by a plateau, i.e., a stage where the residual we monitor at each iteration does not decrease significantly over a certain number of iterations depending on the channel length. 
	 
	 The issue of the plateau has been observed by other authors as well, even for other kinds of Stokes problems (see \cite{olshanskii2010acquired}). In particular, it has been shown that in practice, preconditioning based on the pressure mass matrix clusters most eigenvalues, but not the small outliers. The condition number and solver convergence are visibly improved, but the plateau remains. To tackle this specific issue, the targeted effect of deflation is especially well suited.
	 Effective preconditioners for MINRES can be built by combining ideas from deflation, multigrid and reusing spectral information (see \cite{ramage2021using}), assuming a setting where the leading block of the matrix is ill-conditioned.
	 
	 The same tendency to stagnate (see \Cref{fig:conv1DChL}) appears when we use the CRAIG version of the generalized Golub-Kahan bidiagonalization process \cite{arioli2013generalized}, our solver of choice. This segregated method relies on exploiting the block structure of the indefinite system. In a series of papers \cite{KrDaTaArRu2020,KrSoArTaRu2020,KrSoArTaRu2021}, the authors have shown that this solver exhibits fast convergence when applied to certain saddle point systems and in particular to the Stokes problem on the unit square and a rectangular domain. In an attempt to explain however the plateau for the thin-elongated channel flow problem and thus increasingly slow convergence, we observe that the off-diagonal block $\Am$ becomes increasingly ill-conditioned. In particular, the conditioning varies with certain physical parameters, such as the ratio between the length and width of the channel. A related analysis shows that this ratio impacts the condition number of the Schur complement \cite{chizhonkov2000domain}.

	 The contribution of this article is the definition of a deflation strategy for the CRAIG iterative solver, that extends the framework in \cite{gaul2013framework}. It has been shown that the CRAIG solver is equivalent to the \textit{Conjugate Gradient} (CG) method on the Schur complement equation (see \cite[Chapter~5]{orban2017iterative}). 
	 As such, we draw from the existing knowledge about CG and always consider the effects of our procedures on the Schur complement, either implicitly or explicitly. As pointed out above, we can focus on Stokes problems where the condition number itself may not seem problematic, since this can be improved even with simple preconditioners. Instead, we restrict our attention to the presence of outliers in the spectrum, to better isolate the impact of deflation.
	 
	 Following \cite{gaul2013framework}, we define the operators necessary for deflation and the subsequent correction. 
	 Through projection, a certain subspace can be
	 removed from the range of the matrix. 
	 This subspace is chosen to be the one that the iterative solver has difficulties operating with. 
	 The resulting deflated system is solved quickly, yielding a part of the solution to the initial problem.
	 The complete solution is found through a correction, by which we reintroduce the information initially removed by deflation.
	 The basis for the deflated subspace is typically the set of eigenvectors corresponding to the outliers in the spectrum.
	 The basis itself can be different, so long as the vectors span the required subspace. 
	 As such, the projectors can also be built by other means, such as in \cite{nabben2004comparison}, where  multigrid and domain decomposition approaches are discussed. 
	 
	 Our aim is now to deflate the Schur complement implicitly, without constructing it explicitly.
	 As such, we will not employ its eigendecomposition, but rather the related elliptic singular value decomposition of the off-diagonal block (see \cite{arioli2013generalized}). 
	 Although we define deflation using exact information about the decomposition, in practice it is rarely the case that this is available. 
	 As such, we also study the implications of using approximate vectors. 
	 For computing this information in an iterative way, we extend two already existing algorithms. 
	 The first is the Augmented Lanczos Bidiagonalization by Baglama and Reichel \cite{baglama2005augmented}. 
	 Our extension relies primarily on replacing the original, standard inner-products such that the resulting algorithm realizes the generalized Golub-Kahan bidiagonalization \cite{arioli2013generalized}. Following this change, the output is a set of approximate elliptic singular triplets. For an augmented (not generalized) Golub-Kahan bidiagonalization in the context of the LSQR solver where deflation is carried out using harmonic Ritz values, see \cite{baglama2013augmented}.
	 The second algorithm we extend is inspired by eigCG, which was defined by Stathopoulos and Orginos \cite{stathopoulos2010computing}. 
	 Their version solves a linear system with CG and, at the same time, uses the vectors that the solver builds to compute approximate eigenpairs (also known as Ritz pairs). 
	 With the same objectives, we adapt the CRAIG algorithm of Arioli \cite{arioli2013generalized} such that, as it solves a saddle point system, it also computes approximate elliptic singular triplets.
	 
	 One alternative choice for tackling saddle point systems is to treat the matrix as a single unit and solve the associated system with MINRES \cite{benzi2005numerical}. 
	 It has been shown that CRAIG needs two times less iterations to converge than MINRES, when the latter is paired with a particular block diagonal preconditioner, which makes the two algorithms comparable in terms of cost per iteration \cite{arioli2013generalized}.
	 We perform numerical experiments to test whether this 2:1 ratio still holds if deflation is applied to both algorithms.
	 
	 Throughout the paper, we will illustrate theoretical and practical aspects by making use of two test problems chosen for their simplicity and their specific features. 
	 Both represent a kind of Stokes flow in a long and thin rectangular channel, but with different definitions and discretizations.
	 One of them is designed by us for the purpose of exposing the algorithmic features, while the other is created with the software package Incompressible Flow \& Iterative Solver Software \footnote{http://www.cs.umd.edu/~elman/ifiss3.6/index.html} (IFISS).

	 The paper is structured as follows. In \Cref{sec:defldefi}, we extend the deflation framework presented in \cite{gaul2013framework} to the case of saddle point systems, operating with the individual blocks of the matrix. This deflation is achieved using exact spectral information, following an elliptic singular value decomposition. Next, in \Cref{sec:testPb}, we introduce our test problems, which we use in the remainder of the article for exemplification. We consider some theoretical consequences of deflation using approximate spectral information in \Cref{sec:apxESVD}. \Cref{sec:compESVD,sec:RitzESVD} are devoted to algorithms for computing the elliptic singular value decomposition. In \Cref{sec:minres}, we numerically compare deflated MINRES and deflated CRAIG. We also briefly compare deflation and preconditioning using the pressure mass matrix in \Cref{sec:deflPrec}, before closing with the conclusions in \Cref{sec:conc}.

	 \section{Deflation definition}
	 \label{sec:defldefi}

	 In this section, we introduce the steps to be taken when applying deflation to a saddle point problem, as well as the necessary prerequisites. In addition to the description of the two main steps of projection and correction, we now provide the formal setting and definition. The main tool we will need is the elliptic singular value decomposition (see \cite{arioli2013generalized}). For the time being, we assume the case where this decomposition provides us with exact information, which is subsequently applied via deflation. 
	 
	 Since we require $\Am$ to have full column rank, the system in \Cref{eq:sad_po} has a unique solution. Our solver of choice is CRAIG (given in \Cref{sec:appendix}), which relies on the generalized Golub-Kahan bidiagonalization introduced by Arioli  \cite{arioli2013generalized}. 
	 In particular, we are interested in accelerating its convergence for a class of problems to be described in \Cref{sec:testPb}. 
	 As noted in \cite[Chapter~5]{orban2017iterative}, there is an equivalence between using CRAIG on \Cref{eq:sad_po} and using CG on the Schur complement equation
	 \begin{equation}\label{eq:schur_comp}
	 	\Am ^T \Wm ^{-1} \Am \pv = \Am ^T \Wm ^{-1} \gvv - \rv, 
	 \end{equation}
	 with $ \Sm = \Am ^T \Wm ^{-1} \Am $ being the Schur complement. According to the literature 
	 (see e.g.~\cite{axelsson2014reaching}), CG's convergence depends on the eigenvalue distribution of the matrix it operates on, in this case $\Sm$. 
	 Clustered values are associated with fast convergence, whereas scattered values, particularly the presence of outliers, lead to slow convergence. 
	 Removing these small outliers from the spectrum via deflation can improve the convergence behavior, as seen in \cite{nabben2004comparison}. This represents our core strategy: we employ deflation to achieve a better eigenvalue clustering for the Schur complement, which translates to faster CRAIG convergence. 
	 The challenge is to do this without building $\Sm$ explicitly.
	 
	 One way to achieve this goal is as follows. By performing a Cholesky decomposition of the (1,1)-block, we have $ \Wm = \Lm \Lm ^T$, $ \Wm ^{-1}= \Lm ^{-T} \Lm ^{-1} $ and define 
	 \begin{equation}\label{eq:schur2normal_pieces}
	 	\tilde \Am = \Lm ^{-1} \Am , \quad \yv = \Lm ^T \uv , \quad \fv = \Lm ^{-1} \gvv .
	 \end{equation}
	 We can express the Schur complement equation \eqref{eq:schur_comp} as a normal equation
	 \begin{equation}\label{eq:schur_norm_eq}
	 	\tilde \Am  ^T \tilde \Am  \pv = \tilde \Am  ^T \fv - \rv. 
	 \end{equation}
	 
	 There is a connection between the singular value decomposition of $\tilde \Am$ and the eigendecomposition of $ \Sm$: the right singular vectors of $\tilde \Am$ are the eigenvectors of $\Sm$ and the singular values are the square roots of the eigenvalues. Then, we can implicitly deflate $\Sm$ by relying on $\tilde \Am$. 
	 
	 An alternative to using the SVD of $\tilde \Am$  follows an equivalence given in \cite{arioli2013generalized}. Note that the author considers a more general definition of $\tilde \Am $ as  $\tilde \Am = \Wm ^{-\frac{1}{2}} \Am \Nm ^{-\frac{1}{2}}$, for a symmetric and positive definite matrix $\Nm \in \bR ^{n \times n}$. 
	 
	 Let the following be a partial elliptic (generalized) SVD of $\Am$:
	 \begin{equation}\label{eq:esvd}
	 	\Am \Vm = \Wm \Um \mSigma ,\quad
	 	\Am ^T \Um = \Vm \mSigma ,\quad
	 	\Um ^T \Wm \Um  = \Id _k ,\quad
	 	\Vm ^T \Vm  = \Id _k,
	 \end{equation}
	 with $k < n$ triplets, and $\Id _k$ being the $k \times k$ identity matrix. Additionally, $\Am = \Wm \Um \mSigma \Vm ^T$ only when $k=n$. 
	 The columns of $\Vm \in \bR ^{ n \times k }$ are the right elliptic singular vectors, the columns of $\Um \in \bR ^{ m \times k }$ are the left elliptic singular vectors, and the entries of the diagonal matrix $\mSigma \in \bR ^{ k \times k }$ are the elliptic singular values. 
	 Note that the left vectors are not orthonormal in the Euclidean sense, but with respect to an inner product and norm defined by $\Wm$
	 \begin{align*}
	 	\innprodM{\uv_i}{\uv_j} &= \uv_i ^T \Wm \uv_j  = 0, \quad \forall i \neq j, \\
	 	\normM{\uv} &= \sqrt{\innprodM{\uv_i}{\uv_j}} =1.
	 \end{align*}
	 
	 The equivalence given in \cite{arioli2013generalized} states that the elliptic SVD (ESVD) of $\Am$ is linked to the SVD of $\tilde \Am$:
	 the values and the right vectors of both decompositions coincide. 
	 For the left vectors we have $\Um  = \Lm^{-T} \tilde  \Um $, with $\Um$ being the left elliptic singular vectors of $\Am $ and $\tilde \Um $ being the left singular vectors of ${\tilde \Am }$. 
	 Finally, considering this equivalence with that concerning ${\tilde \Am }$ and $\Sm$, it follows that we can deflate the Schur complement implicitly, by referring to the ESVD of $\Am$.
	 
	 Given \Cref{eq:esvd}, a partial ESVD of $\Am$ with $k$ triplets, we define the approximate left-inverse of $ \Am $ as
	 \begin{equation} \label{eq:deflM}
	 	\Mm = \Vm \mSigma^{-1} \Um ^T   ,
	 \end{equation}
	 the $ \Wm -$orthogonal projector 
	 \begin{equation} \label{eq:deflP}
	 	\Pm = \Id _m - \Am \Mm,  
	 \end{equation}
	 and the (Euclidean) orthogonal projector 
	 \begin{equation} \label{eq:deflQ}
	 	\Qm = \Id _n - \Mm \Am. 
	 \end{equation}
	 Then, we define the deflated system as 
	 \begin{equation}\label{eq:simpler_defl}
	 	\left[
	 	\begin{matrix}
	 		\Wm & \Am \Qm \\
	 		\Qm ^T \Am ^{T} & \mZ 
	 	\end{matrix}
	 	\right] 
	 	\left[
	 	\begin{matrix}
	 		\hat{ \uv } \\
	 		\hat{ \pv } 
	 	\end{matrix}
	 	\right]
	 	=
	 	\left[
	 	\begin{matrix}
	 		\gvv \\
	 		\Qm ^T \rv 
	 	\end{matrix}
	 	\right],
	 \end{equation}
	 whose solution we correct by
	 \begin{equation} \label{eq:simpler_corr}
	 	\pv = \Qm \hat{ \pv } + \Mm  \gvv - \Mm  \Wm \Mm  ^T \rv , \quad \uv = \Pm ^T \hat{ \uv } + \Mm  ^T \rv,
	 \end{equation}
	 to retrieve the solution of the original problem in \Cref{eq:sad_po}.
	 
	 Deflation differs from preconditioning in one important way. Since $rank(\Qm)=n-k$, applying $\Qm$ does not preserve the solution set, i.e. the systems in  \Cref{eq:simpler_defl} and \Cref{eq:sad_po} are not equivalent, making the correction in \Cref{eq:simpler_corr} necessary. The correction terms add information complementary to that found by solving the deflated system.
	 
	 It is also possible to define a deflated system similar to that in \Cref{eq:simpler_defl} by including the action of $\Pm$, but it can be shown that the resulting Schur complement is the same as that of \Cref{eq:simpler_defl}. With $\Pm \in \bR ^{m \times m}, \Qm \in \bR ^{n \times n}$ and $m>n$, it is preferable to avoid applying $\Pm$ whenever possible, due to its size.
	 
	 \begin{theorem}\label{th:equiv_corr}
	 	Applying deflation as in \Cref{eq:simpler_defl} and correction as in \Cref{eq:simpler_corr} leads to the solution of the initial saddle point system \Cref{eq:sad_po}.
	 \end{theorem}
	 
	 \begin{proof}
	 	The proof follows by inserting the corrected solution from \Cref{eq:simpler_corr} into the initial saddle point system in \Cref{eq:sad_po} and making use of the equalities in \Cref{eq:simpler_defl}, the definitions in \Cref{eq:deflM,eq:deflP,eq:deflQ} and the properties of the projectors
	 	\begin{equation}\label{eq:proj_prop}
	 		\Pm \Wm  = \Wm  \Pm ^T; \quad  \Pm \Am = \Am \Qm;  \quad \Qm\Qm =\Qm.
	 	\end{equation}
	 \end{proof}

	 We conclude this section with a note on alternative choices of information to be used in defining deflation. The intuitive one is to leverage information about the spectrum via a decomposition, as we have shown above. Another choice would be using multigrid operators (see \cite{nabben2004comparison}). If we have a left basis $\Um$ ($\Wm$-orthogonal) and a right one $\Vm$ orthogonal) which are not sets of elliptic singular vectors, then we set $\mSigma = \Um ^T \Am \Vm$. If we additionally require that $\Am \Vm = \Wm \Um \mSigma$, then the resulting $\Mm, \Pm$ and $\Qm$ will satisfy the conditions in \Cref{eq:proj_prop}. As such, following \Cref{th:equiv_corr}, deflation and correction will find the solution of the initial saddle point system. The choice of $\Um$ and $\Vm$ is free, as long as their columns span the intended subspaces, i.e. those spanned by the elliptic singular vectors corresponding to the values that we intend to deflate.

	 \section{Deflation effects on two test problems}
	 \label{sec:testPb}
	 
	 In this section, we introduce the 
	 saddle point test cases in order to motivate our choices regarding deflation. 
	 Our focus is on cases that exhibit slow convergence, such as 
	 systems arising from the discretization of incompressible flow problems in elongated domains. It is known that increasing the ratio between the width and length of the domain degrades the conditioning of the Schur complement \cite{chizhonkov2000domain}. We will see that this is due to more small outliers that appear in the spectrum, which also drift away from the main cluster of values. 
	 We explore different kinds of discretization and geometries. We analyze the spectra associated with these problems and see how deflation leads to faster solver convergence.
	 
	 \subsection{A simplified model of flow in a 1D channel}
	 \label{subsec:1DCh}
	 
\begin{figure}
	 	\centering
	 	\includegraphics[width=.71\textwidth]{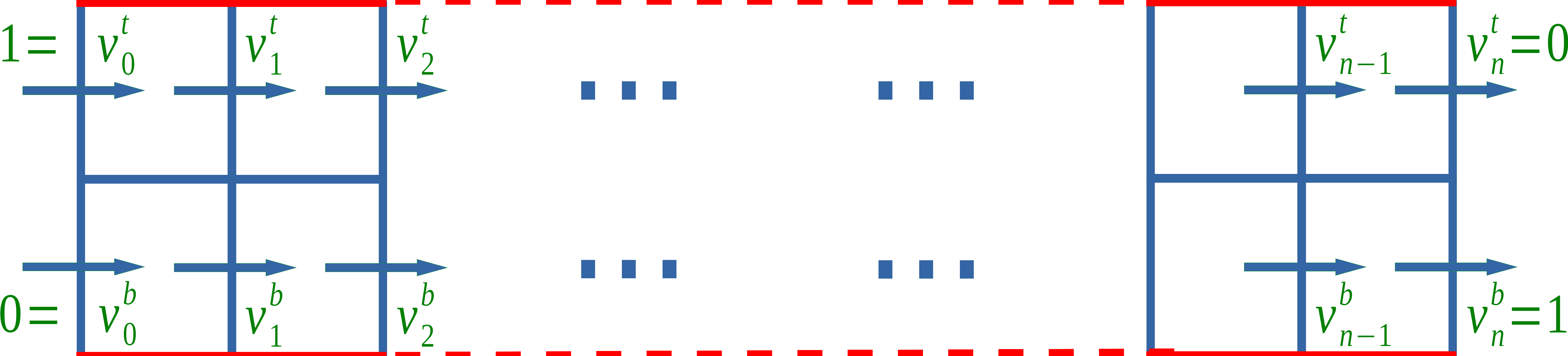}
	 	\caption{1D model of an elongated channel, with inflow only through the upper-left cell ($v ^t _0$) and outflow only through the lower-right cell ($v^b_n$).}
	 	\label{fig:1DChmodel}
	 \end{figure}

	 The Stokes system in its continuous form is given by 
	 \begin{align}
	 	- \Delta \Vec{u}  + \nabla p &= \Vec{f} \label{eq:contStokes}\\
	 	\nabla \cdot  \Vec{u} &=0.
	 \end{align}
	 The first problem is motivated by Stokes flow, but is deliberately constructed to highlight the phenomenon of stagnation in a Krylov method via a simple model system. 
	 To this end, we consider transport in
	 a long and thin channel, as illustrated in \Cref{fig:1DChmodel}.
	 The discrete system results from the classical Marker-and-Cell (MAC) finite difference scheme \cite{harlow1965numerical}.
	 Here the flow enters a channel from the left
	 and leaves at the right.
	 The channel is modeled by a collection of $n$ cells in the horizontal direction, but only two cells vertically.
	 As illustrated in \Cref{fig:1DChmodel}, only the horizontal velocities $v_i^t$ and $v_i^b$, $i=0,1,\ldots ,n$ are used. They represent the velocities at the mid-edge position of the top and bottom cells, respectively.
	 
	 The boundary conditions are set to represent an inflow on the left end.
	 Specifically, we impose an inflow
	 only through the upper half of the boundary, 
	 $v_0^t= 1, v_0^b=0$. 
	 Similarly, outflow is enforced on the right channel end, but only through the lower cell $v_n^t= 0, v_n^b=1$.
	 This setup is chosen to avoid trivial solutions.
	 The velocity values corresponding to these
	 Dirichlet type boundary conditions 
	 can be eliminated and will then
	 appear in the right hand side vector of the reduced linear system.
	 The discrete equations represent a simplified form of momentum conservation taking the form
	 \begin{equation}
	 	\left . \begin{array}{ccc}
	 		- v_{i-1}^t + 4 v_i^t - v_{i+1}^t  - v_i^b & = &  0 \\[1ex]
	 		- v_{i-1}^b + 4 v_i^b - v_{i+1}^b  - v_i^t & = &  0
	 	\end{array} \right \} \text{, for~} i = 1, \ldots, n-1
	 \end{equation}
	 Note that this setup assumes implicitly that the top and bottom walls are characterized by a no-slip boundary condition.
	 Consequently, the resulting velocity block $\Wm$ of the saddle point system in \Cref{eq:sad_po} becomes strictly diagonally dominant and well-conditioned.
	 Using the theorem of Gerschgorin, it is easy to see that its eigenvalues are bounded below by $1$, and its condition number is better than $7$, independent of the number of equations $n$.
	 
	 Next, the lower block row ($[\Am ^T \ \mZ]$ in \Cref{eq:sad_po}) enforces mass conservation for each pair of vertical cells, by requiring
	 \begin{equation}
	 	\begin{array}{ccc}
	 		v_{i}^t  -  v_{i-1}^t + v_{i}^b - v_{i-1}^b & = &  0,
	 	\end{array} \text{~for~} i = 1, \ldots, n
	 \end{equation}
	 The resulting system is as in \Cref{eq:sad_po}.
	 Given this construction, $ \Am $ would be rank-deficient, reflecting the fact that the pressure is determined only up to a constant.
	 To remedy this in the discrete system, we additionally replace the first constraint (on $i=1$) by a convex combination between itself and the last constraint and  
	 the latter is dropped from the system. Once we include the treatment of the boundary conditions and constraints as described above, and we take a row-wise ordering of the unknowns, the entries of the system from \Cref{eq:sad_po} are
	 \begin{equation}\label{eq:entries_sad_po_mat}
	 	\Wm=
	 	\begin{bmatrix*}[r]
	 		4 & -1 &  &  & \vline & -1 &  &  \\
	 		-1 & 4 & -1 &  & \vline &  & -1 &  \\
	 		& -1 & 4 &  \mddots &  \vline &  &  & -1 \\
	 		&  &  \mddots & \mddots & \vline &&&& \mddots \\
	 		\cmidrule(lr) {1-9}
	 		-1 &  &  &  & \vline & 4 & -1 &  & \phantom{\mddots}\\
	 		& -1 & &  & \vline & -1 & 4 & -1  & \phantom{\mddots}\\
	 		& & -1 & &   \vline & & -1 & 4 & \mddots   \\
	 		&&  & \mddots & \vline &&& \mddots & \mddots
	 	\end{bmatrix*}, \quad
	 	\Am= 
	 	\begin{bmatrix*}[r]
	 		0.5 & -1 &    &  \\ 
	 		&  1 & -1 &  \\
	 		&    &  1 &  \mddots  \\
	 		-0.5 &    &    & \mddots \\
	 		\cmidrule(lr) {1-4}
	 		0.5 & -1 &    &  \\ 
	 		&  1 & -1 &  \\
	 		&    &  1 & \mddots   \\
	 		-0.5 &    &    & \mddots 
	 	\end{bmatrix*}, 
	 \end{equation}
	 \begin{equation}
	 	\label{eq:entries_sad_po_vec}
	 	\uv=
	 	\begin{bmatrix*}[c]
	 		v^t _1 \\
	 		\mvdots \\
	 		v^t _{n-1} \\
	 		\cmidrule(lr) {1-1}
	 		v^b _1 \\
	 		\mvdots \\
	 		v^b _{n-1}
	 	\end{bmatrix*},
	 	\quad
	 	\pv=
	 	\begin{bmatrix*}[c]
	 		p _1 \\
	 		\mvdots \\
	 		p _{n-1} 
	 	\end{bmatrix*},
	 	\quad
	 	\gvv= 
	 	\begin{bmatrix*}[c]
	 		v^t _0 =1\\
	 		\mvdots \\
	 		v^t _{n} =0\\
	 		\cmidrule(lr) {1-1}
	 		v^b _0 =0\\
	 		\mvdots \\
	 		v^b _{n} =1
	 	\end{bmatrix*},
	 	\quad
	 	\rv=
	 	\begin{bmatrix*}[c]
	 		1 \\
	 		\mvdots \\
	 		0 
	 	\end{bmatrix*}.
	 \end{equation}
	 
	 In the limit $n \rightarrow \infty$, the channel grows to infinite length, while its width stays constant, or, alternatively, it can be interpreted as staying at constant 
	 length but its width shrinking to zero. 
	 Thus, the linear system is not a discretization of a differential equation on a fixed domain.
	 However, such a change of the domain geometry has an impact on the \textit{inf-sup} condition or equivalently, on the smallest eigenvalue of the Schur complement (see \cite{chizhonkov2000domain}), leading it to approach zero. 
	 
	 Note also that the solution of this particular saddle point system becomes trivial with a change of variables, involving the sum and difference of each vertical pair of velocity variables, as a simple analogue to
	 the stream-function-vorticity formulation of the Stokes problem (see \cite{griebel1998numerical}).
	 
	 Here we do not transform the variables accordingly, but use the original definition of the variables to exhibit and illustrate a typical CRAIG convergence pattern. 
	 When used to solve the initial system from \Cref{eq:sad_po}, a stagnation in error reduction occurs. 
	 This convergence plateau lasts for a number of iterations roughly equal to a quarter of the channel's length, as shown in \Cref{fig:conv1DChL}. 
	 Only when the end of the plateau is reached, fast convergence occurs
	 to then reach the solution quickly.
	 As such, simply by modifying the length of the channel, we can adapt the difficulty of the problem.
	 
	 \subsubsection{Applying deflation on the 1D channel problem}
	 
	 Since the CRAIG algorithm is equivalent to that of the Conjugate Gradient (CG) applied to the Schur complement equation, we can learn more about CRAIG's convergence behavior by analogy with CG. We know from \cite{golub1996matrix} that CG's convergence rate depends on the condition number of the system matrix. For us, this means studying the eigenvalue distribution of the Schur complement or $ \Am$'s elliptic singular values, as defined in \Cref{sec:defldefi}. 
	 
\begin{figure}
	 	\centering	
	 	\begin{subfigure}{\FigWidHalf \textwidth}
    \centering
    \begin{tikzpicture}
    \begin{axis}[
    ymode=log,
    legend pos= south west,
	table/col sep=tab,
    width=\FigWid \textwidth,
	height=\FigHei \textwidth,
    xlabel={CRAIG iterations},
     ylabel={Error norm}
    ]

\addplot+[] table [x=1024X, y=1024Y]{images/csv/conv1DChL.csv}; 
\label{fig:item:conv1DChL_1024} 
	
\addplot+[] table [x=512X, y=512Y]{images/csv/conv1DChL.csv}; 
\label{fig:item:conv1DChL_512} 

\addplot+[] table [x=256X, y=256Y]{images/csv/conv1DChL.csv}; 
\label{fig:item:conv1DChL_256} 

\addplot+[] table [x=128X, y=128Y]{images/csv/conv1DChL.csv}; 
\label{fig:item:conv1DChL_128} 

\end{axis} 
\end{tikzpicture} 
\caption{CRAIG convergence curves. Y-Axis: Energy norm of the relative error (Velocity field).} 
\label{fig:conv1DChL} 
 \label{fig:LongChannelConv}
	 	\end{subfigure}
	 	\begin{subfigure}{\FigWidHalf \textwidth}
    \centering
    \begin{tikzpicture}
    \begin{axis}[
    ymode=log,
    legend pos= south west,
	table/col sep=tab,
    width=\FigWid \textwidth,
	height=\FigHei \textwidth,
    xlabel={Ell. sing. value index},
    ylabel={Magnitude}
    ]
    
    \addplot+[only marks] table [x=L1024X, y=L1024Y]{images/csv/EllSpec1DChL.csv}; 
    \label{fig:item:EllSpec1DChL_L1024} 
    
    \addplot+[only marks] table [x=L512X, y=L512Y]{images/csv/EllSpec1DChL.csv}; 
    \label{fig:item:EllSpec1DChL_L512} 
    
    \addplot+[only marks] table [x=L256X, y=L256Y]{images/csv/EllSpec1DChL.csv}; 
    \label{fig:item:EllSpec1DChL_L256} 

    \addplot+[only marks] table [x=L128X, y=L128Y]{images/csv/EllSpec1DChL.csv}; 
    \label{fig:item:EllSpec1DChL_L128} 
    
    \end{axis} 
    \end{tikzpicture} 
    \caption{Distribution of elliptic singular values.}
    \label{fig:EllSpec1DChL} 
	 	\end{subfigure}	
	 	\caption{Corresponding spectra and CRAIG convergence curves for four 1D channels of lengths 128 \ref{fig:item:conv1DChL_128}, 256 \ref{fig:item:conv1DChL_256}, 512 \ref{fig:item:conv1DChL_512} and 1024 \ref{fig:item:conv1DChL_1024}.}
	 \end{figure}
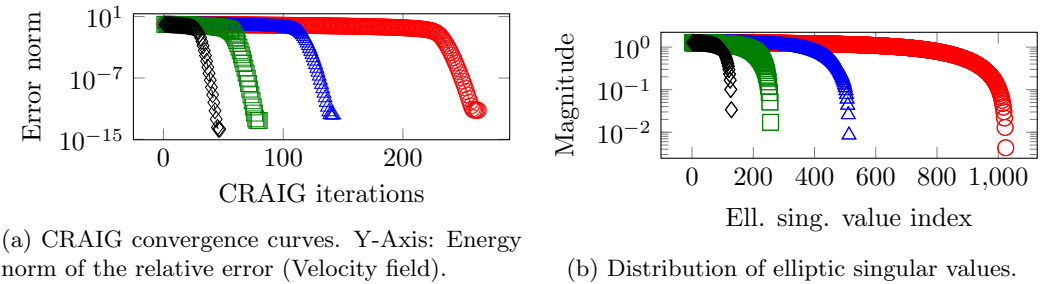

	 In \Cref{fig:EllSpec1DChL}, we illustrate the elliptic spectra of a 1D channel problem as described above, where the channel length varies between 128 and 1024. The largest values are well clustered and away from the origin, but the smallest stray from this cluster. These outliers get more numerous and closer to zero as the channel length increases. In the following considerations, we focus on the case where this length is 512. If the smallest ten values were deflated, then the effective condition number, which uses only the nonzero values, would be smaller and CRAIG would converge in fewer steps. We consider this strategy analogously to the known behavior of CG \cite{nabben2004comparison}.

	 A more detailed look into the problem comes from analyzing the elliptic singular vectors as well, answering the question "Which values and vectors do we really need to deflate?". 
	 It has been shown in \cite{saunders1995solution} that, at each iteration, CRAIG seeks to minimize the error norm of the primal variable $ \uv $ in \Cref{eq:sad_po}, which also holds for the generalized version we use (see \cite{arioli2013generalized}). We consider the initial error $\tilde \ev _1$, i.e. the difference between the true solution and the first iterate computed by CRAIG. By expressing this initial error as a linear combination of the left elliptic singular vectors $\Um$ of $\Am$,
	 \begin{equation}
	 	\label{eq:ini_err}
	 	\Um\zv=\tilde \ev _1, 
	 \end{equation}
	 we can see which of them contribute to this error and how many of them have a large coefficient in $\zv$. These vectors are then good candidates for deflation, since removing them from the problem reduces the error significantly. Deflating vectors which have small (near zero) coefficients in $\zv$ represents an unnecessary expense, as we argue and illustrate below. See \cite{liesen2004convergence} for a similar analysis, but regarding the eigencomponents of the initial error for CG.
	 
	 We return to our 1D channel problem, apply the ideas described above and plot the entries of $\zv$ in \Cref{fig:ESVpresenceCoef1DCh}.
	 Note that in the following examples, we take a descending ordering of the elliptic singular values and of their vectors. As such, the first entries correspond to the most dominant vectors, i.e. those associated with the largest elliptic singular values. The rightmost entries correspond to the vectors associated with the smallest elliptic singular values.
	 
	 \Cref{fig:ESVpresenceCoefFull1DCh} shows that only odd-numbered left elliptic singular vectors are part of the error, while the coefficients of even-numbered ones are close to machine-zero. A black-box eigensolver typically delivers vectors corresponding to consecutive values. As such, if we tried to deflate a group of $ k $ vectors from this problem, only every second of them would actually be necessary.

	 In \Cref{fig:ESVpresenceCoefFilt1DCh}, we only look at the meaningful entries of $\zv$, those above a threshold of $10^{-7}$. It is clear that the highest entries are to the right, which implies that the vectors contributing the most to the error belong to the smallest elliptic singular values. If we choose to deflate these vectors and their corresponding elliptic singular values, CRAIG is certain to converge faster. While it is common in the literature to focus on deflating the vectors associated with the smallest values, we stress that it is important to first analyze whether these vectors have a significant contribution to the error. That is, one could express the initial error or residual as in \Cref{eq:ini_err}, and check whether the corresponding entries in $\zv$ are the largest, as we have seen was the case for the problem considered above (see \Cref{fig:ESVpresenceCoefFilt1DCh}). If the error components mentioned above are not the largest, the plateau will appear at a later moment, after an initial stage of faster convergence. For such problems, deflation may not be
	 necessary, if we only need solutions of limited accuracy, to which we can converge before
	 the onset of the plateau stage.

	 \begin{figure}
	 	\centering	
	 	\begin{subfigure}{\FigWidHalf \textwidth}
    \centering
    \begin{tikzpicture}
    \begin{axis}[
    ymode=log,
    legend pos= south west,
	table/col sep=tab,
	ymin=1e-16,
	xmin=1,
    width=\FigWid \textwidth,
	height=\FigHei \textwidth,
    xlabel={Ell. sing. vector index},
    ylabel={Magnitude}
    ]
\addplot+[only marks, blue] table [x=data1X, y=data1Y]{images/csv/ESVpresenceCoefFull1DCh.csv}; 
\label{fig:item:ESVpresenceCoefFull1DCh_data1} 

\end{axis} 
\end{tikzpicture} 
\caption{All entries of $\zv$.} 
\label{fig:ESVpresenceCoefFull1DCh} 
	 	\end{subfigure}
	 	\begin{subfigure}{\FigWidHalf \textwidth}
    \centering
    \begin{tikzpicture}
    \begin{axis}[
    ymode=log,
    legend pos= south west,
	table/col sep=tab,
	xticklabels={0,1,200,400},
    width=\FigWid \textwidth,
	height=\FigHei \textwidth,
    xlabel={Ell. sing. vector index},
    ylabel={Magnitude}
    ]
	\addplot+[only marks, blue] table [x=data1X, y=data1Y]{images/csv/ESVpresenceCoefFilt1DCh.csv}; 
\label{fig:item:ESVpresenceCoefFilt1DCh_data1} 

\end{axis} 
\end{tikzpicture} 
\caption{ Only entries with magnitudes above $ 10^{-7} $.} \label{fig:ESVpresenceCoefFilt1DCh} 
	 	\end{subfigure}	
	 	\caption{
	 		Entries of the $\zv$ vector (see \Cref{eq:ini_err}), showing how much each left elliptic singular vector contributes to the initial error.1D channel problem (length=512). Vectors listed in descending order of the corresponding elliptic singular values.
	 	}
	 	\label{fig:ESVpresenceCoef1DCh}
	 \end{figure}
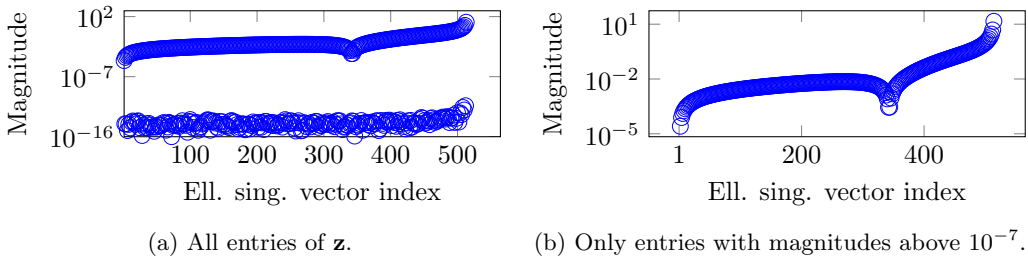
	 
	 \begin{figure}
    \centering
    \begin{tikzpicture}
    \begin{axis}[
    ymode=log,
    legend pos=outer north east,
	table/col sep=tab,
    width=\FigWidSolo \textwidth,
	height=\FigHeiSolo \textwidth,
    xlabel={CRAIG iterations},
     ylabel={Error norm}
    ]
	\addplot+[] table [x=10 smallest X, y=10 smallest Y]{images/csv/conv1DChDefl.csv}; 
\label{fig:item:conv1DChDefl_10 smallest } 
\addlegendentry{10 smallest } 

\addplot+[] table [x=no deflationX, y=no deflationY]{images/csv/conv1DChDefl.csv}; 
\label{fig:item:conv1DChDefl_no deflation} 
\addlegendentry{no deflation} 

\addplot+[] table [x=50 smallestX, y=50 smallestY]{images/csv/conv1DChDefl.csv}; 
\label{fig:item:conv1DChDefl_50 smallest} 
\addlegendentry{50 smallest} 

\addplot+[] table [x=100 smallestX, y=100 smallestY]{images/csv/conv1DChDefl.csv}; 
\label{fig:item:conv1DChDefl_100 smallest} 
\addlegendentry{100 smallest} 

\end{axis} 
\end{tikzpicture} 
\caption{Convergence curves of the CRAIG algorithm for the 1D channel problem after having deflated a number of elliptic singular values of $\Am$. Y-Axis: Energy norm of the relative error (Velocity field).} 
\label{fig:conv1DChDefl} 
\end{figure}  
	 
	 We conclude this subsection with a convergence plot (\Cref{fig:conv1DChDefl}) showing how deflation helps in the case of this 1D channel test case. We compare deflating the smallest 10, 50 and 100 values and take the undeflated problem as reference. Unsurprisingly, deflating more values leads to faster convergence. While it is unrealistic to deflate so many of them, the performance is predictable once we look at the spectrum in \Cref{fig:EllSpec1DChL} and see how the remaining values are distributed. We note that even deflating ten values can be very helpful, especially if we only need a solution of moderate precision.

	 \subsection{IFISS channel problem}
	 \label{subsec:IFISSCh}

	 This test case is generated using the Incompressible Flow \& Iterative Solver Software \footnote{http://www.cs.umd.edu/~elman/ifiss3.6/index.html} (IFISS) package. The particular scenario that we construct is once again that of a long and thin channel domain $\Omega = \left[-1, L \right] \times \left[-1, 1 \right]$, with $L>1$, potentially $L \gg 1$. This time, however, we employ the Finite Elements discretization Q2-Q1. 
	 
	 The steady Stokes problem here is the Poiseuille flow, with the following boundary conditions: a Dirichlet condition on the inflow (left boundary, $\Gamma_{in}= \left\{-1\right\} \times \left[-1, 1 \right]$ ), a characteristic boundary with zero velocity above and below, also known as no-slip condition ($\Gamma_{c}= \left[-1, L \right] \times \left\{-1\right\} \cup \left[-1, L \right] \times \left\{1\right\} $), while the outflow (right, $\Gamma_{out}= \left\{L \right\} \times \left[-1, 1 \right]$) is described by the Neumann condition
	 \begin{equation*}
	 	\begin{split}
	 		\frac{\partial u_x}{\partial x} - p &= 0, \\
	 		\frac{\partial u_y}{\partial x}  &= 0.
	 	\end{split}
	 \end{equation*}
	 For this particular problem definition, we have the analytical solution 
	 \begin{equation}
	 	\label{eq:contStokesSol}
	 	\begin{cases}
	 		u_x = 1- y^2, \\
	 		u_y = 0, \\
	 		p = -2x + \text{constant}.
	 	\end{cases}
	 \end{equation}
	 See \cite{elman2014finite} for a more detailed description of this reference Stokes problem.
	 
	 \subsubsection{Applying deflation on the IFISS channel problem}
	 
\begin{figure}
	 	\centering	
	 	\begin{subfigure}{\FigWidHalf \textwidth}
    \centering
    \begin{tikzpicture}
    \begin{axis}[
    ymode=log,
    legend pos=outer north east,
	table/col sep=tab,
    width=\FigWid \textwidth,
	height=\FigHei \textwidth,
    xlabel={\ac{GKB} iterations},
     ylabel={Error norm}
    ]

\addplot+[] table [x=50X, y=50Y]{images/csv/convIFISSChL.csv}; 
\label{fig:item:convIFISSChL_50} 

\addplot+[] table [x=20X, y=20Y]{images/csv/convIFISSChL.csv}; 
\label{fig:item:convIFISSChL_20} 

\addplot+[] table [x=10X, y=10Y]{images/csv/convIFISSChL.csv}; 
\label{fig:item:convIFISSChL_10} 

\end{axis} 
\end{tikzpicture} 
\caption{ CRAIG convergence curves. Y-Axis: Energy norm of the relative error (Velocity field).} 
\label{fig:convIFISSChL} 
	 	\end{subfigure}
	 	\begin{subfigure}{\FigWidHalf \textwidth}
    \centering
    \begin{tikzpicture}
    \begin{axis}[
    ymode=log,
    xtick={0,500,1000,1800},
    legend pos= south west,
	table/col sep=tab, 
    width=\FigWid \textwidth,
	height=\FigHei \textwidth,
    xlabel={Ell. sing. value index},
    ylabel={Magnitude}
    ]

    \addplot+[only marks] table [x=L50X, y=L50Y]{images/csv/EllSpecIFISSChL.csv}; 
    \label{fig:item:EllSpecIFISSChL_L50} 

	\addplot+[only marks] table [x=L20X, y=L20Y]{images/csv/EllSpecIFISSChL.csv}; 
    \label{fig:item:EllSpecIFISSChL_L20} 
    
	\addplot+[only marks] table [x=L10X, y=L10Y]{images/csv/EllSpecIFISSChL.csv}; 
    \label{fig:item:EllSpecIFISSChL_L10} 

\end{axis} 
\end{tikzpicture} 
\caption{Distribution of elliptic singular values.}
 
\label{fig:EllSpecIFISSChL} 
	 	\end{subfigure}	
	 	\caption{Corresponding spectra and CRAIG convergence curves for three IFISS channels of lengths 10 \ref{fig:item:convIFISSChL_10}, 20 \ref{fig:item:convIFISSChL_20} and 50 \ref{fig:item:convIFISSChL_50}.}
	 \end{figure}
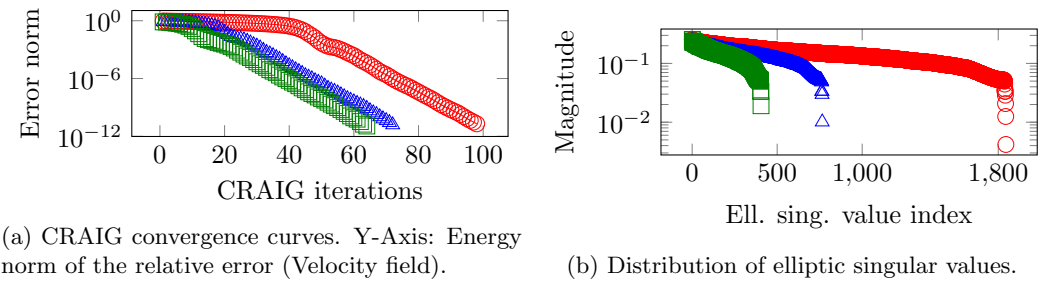

	 As in \Cref{subsec:1DCh}, when solving this problem with CRAIG, we can see in \Cref{fig:convIFISSChL} how the convergence curve starts with a plateau having a length depending on the length of the channel.
	 
	 The elliptic spectrum in \Cref{fig:EllSpecIFISSChL} shows how, for different channel lengths, the main cluster of values is contained within the same interval. Only the smallest values become outliers, getting more numerous and further away from the main cluster. This latter feature explains why the problem gets increasingly more difficult to solve, with a longer plateau. The first feature, concerning the cluster, explains why the convergence curve after the plateau is the same in all three cases. We focus only on the channel of length 20 and will consider the five smallest values as candidates for deflation.

\begin{figure}
	 	\centering	
	 	\begin{subfigure}{\FigWidHalf \textwidth}
    \centering
    \begin{tikzpicture}
    \begin{axis}[
    ymode=log,
    xtick={0, 200, 500, 715 },
    legend pos= south west,
	table/col sep=tab,
	xmin=1,
    width=\FigWid \textwidth,
	height=\FigHei \textwidth,
    xlabel={Ell. sing. vector index},
    ylabel={Magnitude}
    ]
	\addplot+[only marks, blue] table [x=data1X, y=data1Y]{images/csv/ESVpresenceCoefFullIFISSCh20.csv}; 
\label{fig:item:ESVpresenceCoefFullIFISSCh20_data1} 

\end{axis} 
\end{tikzpicture} 
\caption{All entries of $\zv$.} 
\label{fig:ESVpresenceCoefFullIFISSCh20} 
	 	\end{subfigure}
	 	\begin{subfigure}{\FigWidHalf \textwidth}
    \centering
    \begin{tikzpicture}
    \begin{axis}[
    ymode=log,
    ymin=1e-5,
    xtick={1, 200, 500, 715 },
    legend pos= south west,
	table/col sep=tab,
    width=\FigWid \textwidth,
	height=\FigHei \textwidth,
    xlabel={Ell. sing. vector index},
    ylabel={Magnitude}
    ]
	\addplot+[only marks, blue] table [x=data1X, y=data1Y]{images/csv/ESVpresenceCoefFullIFISSCh20.csv}; 
\label{fig:item:ESVpresenceCoefFiltIFISSCh20_data1} 

\end{axis} 
\end{tikzpicture} 
\caption{ Only entries with magnitudes above $ 10^{-7} $.} 
\label{fig:ESVpresenceCoefFiltIFISSCh20} 
	 	\end{subfigure}	
	 	\caption{ Entries of the $\zv$ vector (see \Cref{eq:ini_err}), showing how much each left elliptic singular vector contributes to the initial error. IFISS channel problem (length=20). Vectors listed in descending order of the corresponding elliptic singular values.     }
	 	\label{fig:ESVpresenceCoefIFISSCh20}
	 \end{figure}
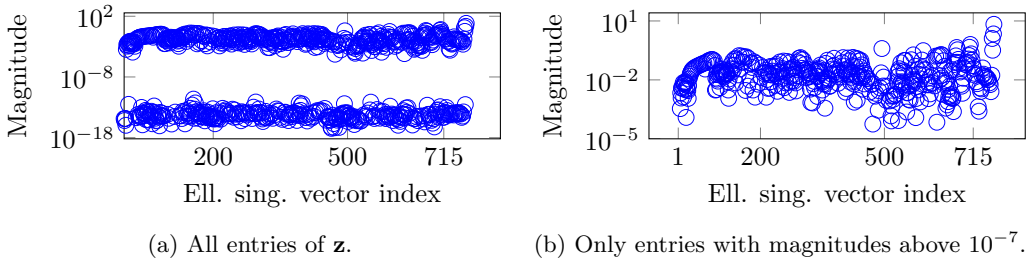

	 We also take a look at the elliptic singular vectors, and how present they are in the initial error for this particular problem. Their coefficients in the linear combination (\Cref{eq:ini_err}) are plotted in \Cref{fig:ESVpresenceCoefIFISSCh20}. We can see in \Cref{fig:ESVpresenceCoefFullIFISSCh20} that many coefficients are close to machine-zero, including that corresponding to the vector associated with the fourth smallest elliptic singular value. 
	 In \Cref{fig:ESVpresenceCoefFiltIFISSCh20}, where we focus only on the entries greater than $10^{-7}$, we notice how the rightmost coefficients are larger than the rest, supporting the idea of deflating the vectors of the smallest elliptic singular values.
	 
\begin{figure}
    \centering
    \begin{tikzpicture}
    \begin{axis}[
    ymode=log,
    legend pos=outer north east,
	table/col sep=tab,
    width=\FigWidSolo \textwidth,
	height=\FigHeiSolo \textwidth,
    xlabel={\ac{GKB} iterations},
    ylabel={Error norm}
    ]
	\addplot+[] table [x=no deflationX, y=no deflationY]{images/csv/convIFSSChDefl.csv}; 
\label{fig:item:convIFSSChDefl_no deflation} 
\addlegendentry{no deflation} 

\addplot+[] table [x=5 smallestX, y=5 smallestY]{images/csv/convIFSSChDefl.csv}; 
\label{fig:item:convIFSSChDefl_5 smallest} 
\addlegendentry{5 smallest} 

\addplot+[] table [x=50 smallestX, y=50 smallestY]{images/csv/convIFSSChDefl.csv}; 
\label{fig:item:convIFSSChDefl_50 smallest} 
\addlegendentry{50 smallest} 

\end{axis} 
\end{tikzpicture} 
\caption{Convergence curves of the CRAIG algorithm for the IFISS channel problem after having deflated a number of elliptic singular values. Y-Axis: Energy norm of the relative error (Velocity field).} 
\label{fig:convIFSSChDefl} 
\end{figure}
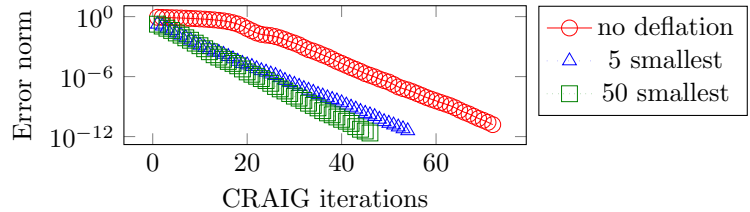

	 Simply from analyzing the spectrum in \Cref{fig:EllSpecIFISSChL} and the coefficients in \Cref{fig:ESVpresenceCoefFiltIFISSCh20}, we can tell that removing the five smallest values should reduce the condition number significantly, accelerating CRAIG's convergence. Deflating more values however, is not expected to bring much additional improvement, since the rest of the values are already well clustered. We test these assumptions and deflate first five, then 50 values, plotting CRAIG's convergence history in \Cref{fig:convIFSSChDefl}. Deflating five values completely removed the plateau. Deflating 50 reduced the number of iterations more than the case of five, but not by much. This is understandable if we consider that in \Cref{fig:ESVpresenceCoefFullIFISSCh20}, more than half of the values after index 715 are already 0, as if they had already been deflated. The rest are nonzeros, but are small nonetheless, not contributing much to the error. Deflation sets them to zero and leads to the additional gains shown in \Cref{fig:convIFSSChDefl}.
	 
	 \section{Deflation using approximate elliptic singular triplets}
	 \label{sec:apxESVD}

	 We have seen how deflation using exact information about the spectrum of the problem can be beneficial by leading to a better spectral clustering and reducing the number of iterations necessary for convergence. Finding such exact information, however, can prove to be prohibitively expensive. It is often the case that one seeks a solution of moderate accuracy and uses an iterative method with a fairly high tolerance parameter. It seems unreasonable then to invest so much in finding exact spectral information, only to accelerate finding an inexact solution to the linear system. Luckily, we can show that this is not necessary and that satisfactory results can be achieved even with approximate knowledge of the deflation vectors. We analyze the spectrum of the deflated matrix to see how it changes when we vary the accuracy of the deflation vectors. Then, we compare the numerical performance of this approximate deflation with that of its exact counterpart.
	 
	 In the literature, the effects of using approximate deflation have been studied for particular definitions of the deflation operator and choice of linear solver. The spectrum of the deflated matrix and its influence on the convergence behavior of GMRES has been analyzed in \cite{sifuentes2013gmres}. In \cite{giraud2006sensitivity}, the authors examine spectral preconditioners, including deflation and coarse grid correction, and how the resulting spectrum changes as a function of perturbation. Subdomain deflation combined with variants of the Incomplete Cholesky decomposition as preconditioner and their effect on the spectrum are discussed in \cite{van2010deflated}.
	 
	 In this paper, we are interested in how approximate spectral information affects the theory presented in \Cref{sec:defldefi}. To understand where the inexactness originates from, let us review the stages of computing the SVD, assuming exact arithmetic. 
	 
	 To find the SVD of a general matrix, three steps are required: a reduction to bidiagonal form, an SVD (diagonalization) of the resulting bidiagonal matrix and a multiplication of the matrices from the two previous steps. One algorithm used for the first step is that of Lanczos \cite{baglama2005augmented}. If we need the singular vectors to be orthogonal in a non-Euclidean way, we need to replace Lanczos's algorithm with the generalized Golub-Kahan bidiagonalization \cite{arioli2013generalized}, which redefines the inner products by using SPD matrices other than the identity. In particular, we require $\Wm$-orthogonal left vectors. With this change, what we find are the elliptic singular triplets \cite{arioli2013generalized}.
	 
	 In the following relations, we will use the subscript $b$ to denote matrices coming from the first stage (the bidiagonalization), and the subscript $s$ for those coming from the second one (the SVD). The final matrices, used to build the deflation operators $\Pm$, $\Qm$ and $\Mm$ have no subscript.
	 
	 We focus on approximate elliptic singular triplets resulting from a decomposition which includes a partial generalized Golub-Kahan bidiagonalization. Assume we stop the bidiagonalization after $\eta < n $ steps. 
	 Then, for $\Am \in \bR ^{m \times n}$, we have $\Um_b \in \bR ^{m \times \eta}$, $\Vm_b \in \bR ^{n \times \eta}$, $\Bm \in \bR ^{\eta \times \eta}$, with
	 \begin{equation}\label{eq:partGK}
	 	\Am \Vm_b = \Wm \Um_b \Bm  ,\quad
	 	\Am^T \Um_b = \Vm_b \Bm^T  + \cE ,\quad
	 	\Um_b^T \Wm \Um_b = \Id_\eta ,\quad
	 	\Vm_b^T \Vm_b = \Id_\eta .
	 \end{equation}
	 We refer to $\cE$ as the residual term, which vanishes only when $\eta = n$ (and in exact arithmetic), leading to $\Am = \Wm \Um_b \Bm \Vm_b^T$. More importantly, when $\eta = n$ the spectrum of $\Bm$ is the same as that of $\Am$. However, if $\eta <n $, then it is only an approximation which improves as $\eta \rightarrow n$ and/or if the bidiagonalization includes a restarting mechanism. Following \cite{baglama2005augmented}, we can further describe this term as 
	 \begin{equation}\label{eq:partGKerr}
	 	\cE=\rv _\eta \ev _\eta^T,
	 \end{equation}
	 where $\rv _\eta$ is the residual vector used to continue the bidiagonalization: $ \beta_{\eta+1}=\normEu{\rv _\eta}$, $\vvv_{\eta+1}=\rv _\eta / \beta_{\eta+1}$ and $\ev _\eta \in \bR ^ \eta, \ev _\eta  = \{ 0, ..., 0, 1 \} ^ T$; i.e., $\ev _\eta$ is taken to be the last vector of the standard basis for $\bR ^ {\eta \times \eta }$.
	 
	 As stated before, to get our approximate elliptic singular triplets, we perform the SVD of $\Bm$, keep only the $k$ $target$ triplets
	 \begin{equation}\label{eq:partSVD}
	 	\Bm \Vm_s= \Um_s \mSigma  ,\quad
	 	\Vm_s^T \Vm_s = \Id_k   ,\quad
	 	\Um_s^T \Um_s = \Id_k.
	 \end{equation}
	 Above, $target$ denotes which triplets we aim to deflate: either the first (largest elliptic singular values) or the last (smallest).
	 
	 The SVD matrices are then multiplied with those from the partial bidiagonalization
	 \begin{equation}
	 	\Um = \Um_b \Um_s  ,\quad
	 	\Vm = \Vm_b \Vm_s.
	 \end{equation}
	 
	 The approximate elliptic SVD (ESVD) is then
	 \begin{equation}\label{eq:apxESVD}
	 	\Am \Vm = \Wm \Um \mSigma  ,\quad
	 	\Am^T \Um = \Vm \mSigma  + \cE \Um_s  ,\quad
	 	\Vm^T \Vm = \Id_k   ,\quad
	 	\Um^T \Wm \Um = \Id_k.
	 \end{equation}
	 
	 As a side note, we remark that it is also possible to have an error term $\Em$ such that $\Am \Vm = \Wm \Um \mSigma +\Em$, if we consider a different kind of partial bidiagonalization, or simply other sources of error. While they are not considered in this present work, such cases represent interesting directions for further research. 
	 
	 For the residual term in \Cref{eq:apxESVD}, we obtain
	 \begin{align}
	 	\cE  \Um_s = \beta_{\eta+1} \vvv_{\eta+1} \ev _\eta^T \Um_s.
	 \end{align}
	 It is known that the last components of the eigenvectors of the Lanczos tridiagonal matrix can be used to track progress when approximating eigenpairs of a symmetric, positive definite matrix (see \cite{golub2000large,paige1980accuracy}). Given the link with the bidiagonalization, this observation can be used for the SVD as well, as done in \cite{baglama2005augmented}, leading to a stopping criterion for the iterative procedure needed to compute such approximate triplets. The GK bidiagonalization we use here is a generalization, but as seen above, it is still true that the residual term depends on the last row of $\Um_s$. With more restarts and/or a larger subspace generated by the GK bidiagonalization, this residual term diminishes.
	 
	 We have seen how a partial bidiagonalization leads to a residual term in the ESVD. We turn now to the deflated matrix $\Pm\Am$, where $\Pm= \Id_m - \Wm \Um \Um^T$, to see how it changes when relying on approximate spectral information. 
	 \begin{align}\label{eq:apxDefl}
	 	\Pm\Am &= \Am - \Wm \Um \Um^T \Am = \Am - \Wm \Um ( \mSigma \Vm^T + \Um_s^T \cE^T ) \\
	 	&= \Am - \Wm \Um \mSigma \Vm^T - \Wm \Um \Um_s^T \cE^T  
	 \end{align}
	 The partial bidiagonalization also affects the term $\Wm \Um \mSigma \Vm^T$. If the deflation vectors were exact, we would shift $k$ $target$ values to 0, while leaving the rest as they were. Neither effect can be counted on in the presence of inexactness, here due to approximating the elliptic singular triplets of $\Am$ from the small bidiagonal matrix $\Bm$. The spectrum will include $k$ zeros, but not necessarily from deflating the intended values. Also, the other values may suffer perturbations. Both effects diminish as the accuracy of the deflation vectors increases. If $\eta=n$ and we assume exact arithmetic, we do indeed get the intended triplets. In general, neither condition holds and we only get an approximation that improves as $\eta \rightarrow n$, or via a restarted algorithm. These changes to the spectrum then propagate, leading linear solvers to converge differently as we will explore numerically in \Cref{sec:compESVD}.

	\section{Computing the elliptic singular triplets}
	\label{sec:compESVD}

	In this section, we describe the algorithms by which one can find an approximate elliptic SVD (ESVD) and provide numerical illustration of the points brought up in \Cref{sec:apxESVD} regarding the influence of the eigensolver on the linear solver with respect to accuracy. 
	The first approach to computing the ESVD (\Cref{alg:ESVD1}) uses already existing tools while the second approach (\Cref{alg:ESVD2}) relies on a reformulation of the generalized Golub-Kahan bidiagonalization algorithm, not as part of a linear solver, but as part of an eigensolver.
	
	In principle, the ESVD of $ \Am $ differs from a regular SVD only in the requirement that the left singular vectors should be $ \Wm$-orthogonal, instead of Euclidean. As such, we can "integrate" this target space of orthogonality into the problem, and then compute the SVD. According to the description in \Cref{sec:defldefi}, we have the following \Cref{alg:ESVD1}.
	
	\begin{algorithm} 
		\caption{Elliptic Singular Value Decomposition 1}
		\label{alg:ESVD1}
		\begin{algorithmic}[1] 	
			\REQUIRE \hspace{-0.2cm}: $ \Lm $ - Cholesky factor of $ \Wm$, lower triangular matrix; 
			$ \Am $ - rectangular matrix.
			\ENSURE \hspace{-0.2cm}: $ \Um $ - left elliptic singular vectors; 		$ \mSigma $ - elliptic singular values; 	$ \Vm $ - right elliptic singular.
			\STATE Compute a (partial) SVD of $ \Lm^{-1}\Am $: get left singular vectors $ \tilde \Um $, singular values $ \mSigma $ and right singular vectors $ \Vm $.
			\STATE Transform:  $ \Um = \Lm^{-T} \tilde  \Um $	
		\end{algorithmic}
	\end{algorithm}
	
	\noindent
	This algorithm relies on the ability to factorize $ \Wm.$ If that is the case, it enables us to find the ESVD, and also reuse these factors when we call CRAIG to solve the linear system. If the Cholesky decomposition cannot be performed, \Cref{alg:GKB_ESVD} is a viable alternative. 
	
	The second approach of finding a partial, approximate ESVD follows the stages described in \Cref{sec:apxESVD}. It is a nested iteration, where the outer algorithm (\Cref{alg:ESVD2}) carries out a partial bidiagonalization (\Cref{alg:GKB_ESVD}) at every step. As seen in \Cref{sec:apxESVD}, the accuracy of the ESVD depends to a significant extent on the bidiagonalization, i.e. on how well the spectrum of $\Bm$ approximates that of $\Am$. Due to memory limitations, the bidiagonalization is not performed only once with $n$ steps, but rather several times, with $\eta<n$ steps, in a restarted fashion. The key here is that with each restart, the bidiagonalization improves the results of the previous one via an augmentation mechanism, which allows it to explore a new portion of the search space.
	
	To define our algorithms, we extend the developments presented in  \cite{baglama2005augmented}. 
	The strategy presented there is to augment the bidiagonalization search space at each restart with the current approximation of the $ k $ elliptic singular triplets. From the second call on, the bidiagonalization generates $ \eta-k-1$ triplets, from which information is extracted by a subsequent SVD to refine the existing $ k $ elliptic singular triplets. The  presentation in \cite{baglama2005augmented} uses the Lanczos bidiagonalization, while we need the generalized GK algorithm, in order to have $ \Wm $- orthogonal left vectors.

	\begin{algorithm}
		\caption{Elliptic Singular Value Decomposition 2}
		\label{alg:ESVD2}
		\begin{algorithmic}[1] 	
			\setcounter{ALC@unique}{0}
			\REQUIRE  \hspace{-0.2cm}: 
			$ \Am $ - rectangular matrix, 
			$ \Wm $ - SPD matrix,
			$ \vvv $ - starting vector, 	
			$ \eta $ - maximum search subspace dimension, 
			$ k $ - number of desired triplets, 		
			$ target $ - choose between $ smallest $ or $ largest $ triplets, 
			$ tol $ - tolerance, 
			$ max\_iter $ - maximum number of iterations.	
			\ENSURE   \hspace{-0.2cm}: 
			$ \Um $ - approximate left elliptic singular vectors,  
			$ \mSigma $ - approximate elliptic singular values,  
			$ \Vm $ - approximate elliptic singular right vectors. 
			\vspace{.5cm}
			
			\STATE 	Call \Cref{alg:GKB_ESVD} with a starting vector $ \vvv $, $\uv _0 = \mZ $, and passing $ \Am, \Wm, $ and $ \eta $. Get $ \Um_b, \Bm, \Vm_b $ and residual $ \rv $. \label{alg_line:1stGKBcall}
			
			\FOR	    { $ i = 1 : max\_iter $  }		
			\STATE Compute the SVD of $ \Bm  $ for the $ k $ {\it target} triplets; store them in  $ \Um_s, \mSigma _i $ and $ \Vm_s $.
			\STATE  $ \Um _i= \Um_b \Um_s $;
			$\quad$  $ \Vm _i= \Vm_b \Vm_s $;
			$\quad$ $  \beta = \normEu{\rv}$;
			\STATE $\hat{\uv} =$ the last row of $\Um_s$; 
				$\quad$ $ \mu =\max(| \hat{\uv} |);$ 
			$\quad$ $\sigma=$ the first entry of $\mSigma _i$;
			\IF	   { $  (\beta  \mu) / \sigma  < tol $ OR $  i  = max\_iter $}
			\STATE \textbf{break}
			\ENDIF	
			
			\STATE  $ \vvv =\rv / \beta  $;
			$\quad$  $ \Vm_i = [\Vm_i \  \vvv] $;
			$\quad$  $ \uv = \Am\vvv  $;
			$\quad$  $ \dv = \Um^T _i \uv  $;		
			$\quad$  $ \uv = \Wm^{-1}  \uv - \Um _i \dv $; \label{alg_line:only_ortho}
			$\quad$  $ \mSigma _i= [\mSigma _i \ \dv]$
			$\quad$  $ \alpha = \normM{\uv}$;
			$\quad$  $ \mSigma _i= 
			\begin{bsmallmatrix}
				\mSigma _i \\
				\alpha \ev_{k+1}^T
			\end{bsmallmatrix}
			$;
			$\quad$ $ \uv= \uv/ \alpha$;
			$\quad$ $ \Um _i=[\Um _i \ \uv]$.
			\STATE Compute residual $ \vvv = \Am^T \uv - \alpha \vvv $;
			\STATE Reorthogonalization:	$ \vvv= \vvv - \Vm_i \Vm_i^T$ $\vvv $;
			\STATE $  \beta = \normEu{\vvv}$;
			$\quad$ $ \uv _0 =$ the last column of $\Um _i$;
			\STATE 	Call  \Cref{alg:GKB_ESVD} with starting vector $ \vvv $,  $ \eta - k - 1$ steps and passing $ \Am,  \Wm $ and $\uv _0 $. Optionally, pass $ \Um _i$ and $ \Vm _i$ for additional reorthogonalization. Get $ \Um_b, \Bm, \Vm_b $ and residual $ \rv $.
			\STATE  Extend $\Um _b=[\Um _i\ \Um _b], $  $\Vm _b=[\Vm _i \ \Vm _b] $ and $\Bm =  \begin{bsmallmatrix}
				\mSigma _i & \mZ \\
				\mZ  & \Bm 
			\end{bsmallmatrix}$. Set $ \Bm _{k+1, k+2}  = \beta$.
			\ENDFOR
			\STATE $\Um = \Um _{i}; \quad \Vm = \Vm _{i}; \quad \mSigma = \mSigma _{i}; $
		\end{algorithmic}
	\end{algorithm}

	\begin{algorithm}
		\caption{ Restarted Augmented Golub-Kahan Bidiagonalization }
		\label{alg:GKB_ESVD}
		\begin{algorithmic}[1] 	
			\setcounter{ALC@unique}{0}
			\REQUIRE \hspace{-0.2cm}:  $ \Am $ - rectangular matrix; 
			$ \Wm $ - SPD matrix;
			$ \vvv $ - right starting vector;
			$ \uv _0 $ - left starting vector;
			$ \eta $ - number of steps.\\
			Optional, for additional reorthogonalization: 
			$ \Um $ -  approximate left elliptic singular vectors;
			$ \Vm $ - approximate right elliptic singular vectors. \\
			\ENSURE \hspace{-0.2cm}:   $ \Um_b $ - left vectors ($\Wm$-orthogonal);
			$ \Bm $ - upper bidiagonal matrix;
			$ \Vm_b $ - right vectors (orthogonal);
			$ \rv $ - residual vector
			\vspace{.5cm}
			
			\STATE $  \beta _1= \normEu{\vvv}$; 	$\quad \vvv _1 = \vvv/\beta _1$; $\quad \Vm _1 = \vvv _1$; $\quad \uv= \Wm^{-1}\Am\vvv _1$;
			\IF	    { $ \uv _0 \neq \mZ $   }   
			\STATE $ \Wm $-orthogonalization:	$ \uv= \uv - \beta _1 \uv _0$; \label{alg_line:augmGKB}
			\ENDIF
			\IF	    { $ \Um  $ is given }   
			\STATE $ \Wm $-reorthogonalization:	$ \uv= \uv - \Um \Um^T \Wm \uv $;   
			\ENDIF
			\STATE $  \alpha _1= \normM{\uv}$; $\quad \Bm _1 = \alpha _1$; 	$\quad \uv _1 = \uv/ \alpha _1$; $\quad \Um _1 = \uv _1$;
			\FOR	    { $ j = 1 : \eta $   }
			\STATE $ \vvv =\Am^T \uv _j - \alpha _j \vvv _j;  $
			\STATE Reorthogonalization:	$ \vvv= \vvv - \Vm_{j} \Vm_{j}^T \vvv $;
			\IF	    { $ \Vm $  is given }   
			\STATE Reorthogonalization:	$ \vvv= \vvv - \Vm   \Vm ^T \vvv $;  
			\ENDIF
			\IF	    { $ j < \eta $  }
			\STATE $  \beta _{j+1}= \normEu{\vvv}$;	$\quad \vvv _{j+1}= \vvv/ \beta _{j+1}$; $\quad \Vm _{j+1}= [\Vm _{j} \ \vvv _{j+1}] $;
			\STATE $\Bm _{j+1} = [\Bm _{j} \ \beta _{j+1} \ev_j ]$
			$\quad \uv= \Wm^{-1}\Am\vvv _{j+1}- \beta _{j+1} \uv _{j}$;
			\STATE $ \Wm $-reorthogonalization:	$ \uv= \uv - \Um_{j} \Um_{j}^T \Wm \uv $;
			\IF	    { $ \Um  $ is given }   
			\STATE $ \Wm $-reorthogonalization:	$ \uv= \uv - \Um \Um^T \Wm \uv $;   
			\ENDIF
			\STATE $  \alpha _{j+1}= \normM{\uv}$ ; 
			$\quad \uv _{j+1}= \uv/ \alpha _{j+1}$;
			$\quad \Um _{j+1}  =[ \Um _{j}  \ \uv _{j+1}]$;
			$ \quad$
			$ \Bm _{j+1}=$ 
			$\begin{bsmallmatrix}
				\Bm _{j+1} \\
				\alpha _{j+1} \ev_{j+1}^T
			\end{bsmallmatrix};$
			\ENDIF
			\ENDFOR
			\STATE 	$ \rv= \vvv; \quad \Vm _b = \Vm _{j}; \quad \Um _b = \Um _{j}; \quad  \Bm  = \Bm _{j};$ 
		\end{algorithmic}
	\end{algorithm}

	We note that, owing to the symmetry of $\Wm$ and $\Id$, it is possible for the bidiagonalization to build sets of vectors that are $ \Wm$- and Euclidean orthogonal using short recurrences. As such, after each restart, only the first computed vector $\uv$  needs to be explicitly orthogonalized against all those in $\Um$ (\Cref{alg_line:only_ortho}). The subsequent ones then reference only the immediate predecessor, while staying orthogonal to all the previous ones. 
	
	The bidiagonalization is augmented in \Cref{alg_line:augmGKB}. The first time, \Cref{alg:GKB_ESVD} is called with $\uv _0 = \mZ $ (\Cref{alg_line:1stGKBcall}), since there is no information to augment the bidiagonalization with. 
	Subsequent calls also send the last column of the matrix $\Um$ (the current approximate left singular vectors), as $\uv _0$ which implicitly represents the augmentation. Thanks to the short recurrences mentioned above, we do not need the entire matrix, and even the last column is only referenced once. 
	This way, the augmentation in \Cref{alg:GKB_ESVD} has only a negligible cost, while ensuring that each call finds new information used to refine $\Um, \mSigma$ and $\Vm$ in the outer iteration (\Cref{alg:ESVD2}).

	The stopping condition in \Cref{alg:ESVD2}  depends on the term $\cE \Um_s $ from \eqref{eq:apxESVD}, computed in a simpler way, based on scalars, as in \cite{baglama2005augmented}. The quantity $\beta \mu \rightarrow 0 $ as the iterations progress and $\sigma$ is used in order to achieve a relative, problem-specific value, since it approximates the largest elliptic singular value. This corresponds to what in \cite{venkovic2020comparative} is referred to as eigenresidual. Unfortunately, our tests, like those in \cite{stathopoulos2010computing} and \cite{venkovic2020comparative}, show that the eigenresidual's Euclidean norm does not decrease monotonically, but rather in an irregular way. This happens predominantly in the beginning of the iterative process, stabilizing afterwards. Such irregularities would make us think that the approximations do not improve in quality every iteration, making it difficult to know when to stop the algorithm if we only need rough approximations. However, the theory about the convergence of the restarted Lanczos process in \cite{aishima2015global} shows that the error decreases monotonically. Our tests, as well as those in \cite{stathopoulos2010computing} and \cite{venkovic2020comparative} support this finding. The implication is that the eigenresidual's Euclidean norm may not be the best stopping criterion for \Cref{alg:ESVD2} and that it may not accurately reflect the quality of the approximations.

	The reorthogonalization steps in \Cref{alg:GKB_ESVD} and \Cref{alg:ESVD2} could be necessary in practice, as a result of finite-precision computations. In \Cref{alg:GKB_ESVD}, we can reorthogonalize with respect to local vectors ($\Um _{j}$ and $\Vm _{j}$) computed during that call and/or those computed previously, representing the current approximations for the elliptic singular vectors ($\Um$ and $\Vm$). How this is carried out exactly is a problem dependent choice, but the authors of \cite{baglama2005augmented} and \cite{barlow2013reorthogonalization} provide useful suggestions regarding partial or full reorthogonalization. We can choose to work with $\Um _{j}$ and/or $\Vm _{j}$ or a subset of their columns, in order to fine-tune the relationship between them. An important observation is that, when targeting the smallest values, the computations are more sensitive to round-off errors, so full reorthogonalization is the safest choice, in order to maximize the accuracy of the results. For simplicity, we used the Gram-Schmidt procedure, but more robust choices are also possible. As we will see in the numerical tests, deflation does not need exact information in order to be effective. Thus, if we compute our deflation vectors via \Cref{alg:ESVD2}, we may skip some or all of the reorthogonalization steps, provided that the problem is not very ill-conditioned. See \cite{giraud2005loss} for a study comparing the loss of orthogonality between vectors computed by several variants of the Gram-Schmidt method, as well as \cite{rozlovznik2012numerical}, where this behavior is considered for non-Euclidean orthogonality.

	Next, we analyze one of the test cases given in \Cref{sec:testPb}, to see how the spectrum looks like after deflation with approximate vectors and how this changes the way our linear solver, CRAIG, converges. We use the problem of the 1D channel to illustrate the following points. We run \Cref{alg:ESVD2} for 10, 26 and 30 iterations and use the resulting triplets (the 10 smallest) for deflation. The error and perturbation ranges are given in \Cref{tab:ESVD2err}. Let $\mathcal{I}$ be the set of indices denoting the values we want to deflate and $\overline{\mathcal{I}}$ its complement. We use the superscript $*$ to refer to exact triplets. Then, the  error for the deflated vectors is computed as $\normEu{\vvv_i^* - \Vm \Vm^T \vvv_i^*}, i \in \mathcal{I}$, the relative error for the deflated values is $ |(\sigma_i-\sigma^*_i)/ \sigma ^*_i|, i \in \mathcal{I}$, and the relative perturbation of the undeflated values is $ |(\sigma _i-\sigma ^*_i)/ \sigma ^*_i|, i \in \overline{\mathcal{I}} $.
	
	\begin{table}
		\centering
		\caption{Errors of the approximated elliptic singular triplets and perturbations of the spectrum after deflation with triplets from \Cref{alg:ESVD2}.}
		\label{tab:ESVD2err}
		\begin{tabular}{cccc}
			\begin{tabular}[c]{@{}c@{}}ESVD \\ iterations\end{tabular} & \begin{tabular}[c]{@{}c@{}}median and max error\\  (targeted vectors)\end{tabular} & \begin{tabular}[c]{@{}c@{}}median and max relative\\  error (targeted values)\end{tabular} & \begin{tabular}[c]{@{}c@{}}median and max relative\\ perturbation (other values)\end{tabular} \\
			10                                                         & {[}5E-1, 1{]}                                                                   & {[}1, 2{]}                                                               & {[}1E-5, 4E-2{]}                                                                      \\
			26                                                         & {[}2E-4, 1E-2{]}                                                                & {[}1E-6, 7E-2{]}                                                               & {[}4E-12, 4E-2{]}                                                                     \\
			30                                                         & {[}2E-8, 6E-7{]}                                                                & {[}1E-15, 1E-10{]}                                                             & {[}4E-16, 8E-11{]}                                                                   
		\end{tabular}%
	\end{table}
	
\begin{figure}
		\centering	
		\begin{subfigure}{\FigWidHalf \textwidth}
    \centering
    \begin{tikzpicture}
    \begin{axis}[
    ymode=log,
    legend pos= south west,
	table/col sep=tab,
    width=\FigWid \textwidth,
	height=\FigHei \textwidth,
    xlabel={CRAIG iterations},
     ylabel={Error norm}
    ]
	\addplot+[] table [x=without deflationX, y=without deflationY]{images/csv/1DChConvDeflInExact.csv}; 
\label{fig:item:1DChConvDeflInExact_without deflation} 

\addplot+[] table [x=inexact deflation (10 ESVD iterations)X, y=inexact deflation (10 ESVD iterations)Y]{images/csv/1DChConvDeflInExact.csv}; 
\label{fig:item:1DChConvDeflInExact_inexact deflation (10 ESVD iterations)} 

\addplot+[] table [x=inexact deflation (26 ESVD iterations)X, y=inexact deflation (26 ESVD iterations)Y]{images/csv/1DChConvDeflInExact.csv}; 
\label{fig:item:1DChConvDeflInExact_inexact deflation (26 ESVD iterations)} 

\addplot+[] table [x=inexact deflation (30 ESVD iterations)X, y=inexact deflation (30 ESVD iterations)Y]{images/csv/1DChConvDeflInExact.csv}; 
\label{fig:item:1DChConvDeflInExact_inexact deflation (30 ESVD iterations)} 

\addplot+[] table [x=exact deflationX, y=exact deflationY]{images/csv/1DChConvDeflInExact.csv}; 
\label{fig:item:1DChConvDeflInExact_exact deflation} 

\end{axis} 
\end{tikzpicture} 
\caption{ CRAIG convergence curves. The Y-axis tracks the energy norm of the relative error for the velocity field. } 
\label{fig:1DChConvDeflInExact} 
		\end{subfigure}
		\begin{subfigure}{\FigWidHalf \textwidth}
    \centering
    \begin{tikzpicture}
    \begin{axis}[
    ymode=log,
    legend pos= south west,
	table/col sep=tab,
    width=\FigWid \textwidth,
	height=\FigHei \textwidth,
    xlabel={Ell. sing. values index},
    ylabel={Magnitude}
    ]
	\addplot+[] table [x=original spectrumX, y=original spectrumY]{images/csv/1DChSpecDeflInExactZoom.csv}; 
\label{fig:item:1DChSpecDeflInExactZoom_original spectrum} 

\addplot+[] table [x=inexact deflation (10 ESVD iterations)X, y=inexact deflation (10 ESVD iterations)Y]{images/csv/1DChSpecDeflInExactZoom.csv}; 
\label{fig:item:1DChSpecDeflInExactZoom_inexact deflation (10 ESVD iterations)} 

\addplot+[] table [x=inexact deflation (26 ESVD iterations)X, y=inexact deflation (26 ESVD iterations)Y]{images/csv/1DChSpecDeflInExactZoom.csv}; 
\label{fig:item:1DChSpecDeflInExactZoom_inexact deflation (26 ESVD iterations)} 

\addplot+[] table [x=inexact deflation (30 ESVD iterations)X, y=inexact deflation (30 ESVD iterations)Y]{images/csv/1DChSpecDeflInExactZoom.csv}; 
\label{fig:item:1DChSpecDeflInExactZoom_inexact deflation (30 ESVD iterations)} 

\addplot+[] table [x=after exact deflationX, y=after exact deflationY]{images/csv/1DChSpecDeflInExactZoom.csv}; 
\label{fig:item:1DChSpecDeflInExactZoom_after exact deflation} 

\end{axis} 
\end{tikzpicture} 
\caption{Elliptic singular values sorted in descending order. Zoom in on the smallest 40. }
\label{fig:1DChSpecDeflInExactZoom} 
		\end{subfigure}	
		\caption{ Comparison of deflation effects on the 1D channel problem. We consider the case without deflation \ref{fig:item:1DChConvDeflInExact_without deflation}, three cases with inexact deflation using vectors found after a given number of ESVD iterations (10 \ref{fig:item:1DChConvDeflInExact_inexact deflation (10 ESVD iterations)}, 26 \ref{fig:item:1DChConvDeflInExact_inexact deflation (26 ESVD iterations)} and 30 \ref{fig:item:1DChConvDeflInExact_inexact deflation (30 ESVD iterations)}) and the case with exact deflation \ref{fig:item:1DChConvDeflInExact_exact deflation}. }
	\end{figure}
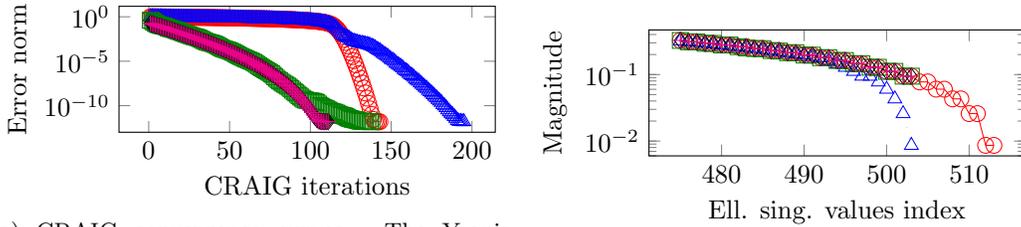

	By construction, the deflated matrix will have a spectrum containing $k=10$ zeros. Depending on the quality of the triplets, these zeros come either from deflating our targeted values or others. By zooming in on the portion of the spectrum that we wanted to deflate (see \Cref{fig:1DChSpecDeflInExactZoom}), we get a detailed view. In the least precise case (10 iterations), we have deflated only some of our targeted values, and since the smallest one is still in place, we cannot expect a smaller condition number or a faster CRAIG convergence. Additionally, the values include a significant perturbation (see \Cref{tab:ESVD2err}). These two effects explain why CRAIG actually needs more iterations in \Cref{fig:1DChConvDeflInExact}. For the case with 26 iterations, all the targeted values get deflated, but there are also a few highly perturbed ones among the rest. The corresponding convergence curve in \Cref{fig:1DChConvDeflInExact} resembles that of exact deflation, until the perturbed values come into play and lead it to deviate. If we run \Cref{alg:ESVD2} for 30 iterations, the resulting information is precise enough for our problem. That is to say, there is no difference between exact and approximate deflation, at least when comparing CRAIG's convergence curves. As for the case with 26 iterations, all the targeted values were successfully deflated, but additionally, the perturbation on the rest approaches machine zero (see \Cref{tab:ESVD2err}).

	How precise should the deflation vectors be? The answer is problem dependent, which means we should think about the matrix spectrum and the coefficients in \Cref{eq:ini_err} related to the initial error. For the 1D channel problem, we have seen in \Cref{sec:testPb} that only every second elliptic singular vector is part of the error. However, the smallest elliptic singular values occur in pairs with (almost) equal values. Then, \Cref{alg:ESVD2} needs some time to correctly identify the vectors, otherwise it may deliver approximations to those who are not part of our problem, while perturbing those that are. An alternative would be to adapt the block Lanczos algorithm (see  \cite{golub1977block}), which is better suited for problems where the values are not simple. Our second test case, the IFISS channel, has a spectrum with outliers that are better separated from the main cluster and from each other (see \Cref{fig:EllSpecIFISSChL}). There, we have found that approximate deflation is as effective as the exact one even with less precise vectors.

	 \section{Approximating elliptic singular triplets from a previous CRAIG solve}
	 \label{sec:RitzESVD}
	 
	 Obtaining approximate information about the spectrum can be done in several ways. One approach is that of finding this information prior to solving the linear system. In this case, we can first use the algorithms given in \Cref{sec:compESVD} to find the deflation vectors, then proceed to solving the deflated linear system. Another approach becomes available if one considers the case of solving several linear systems, where the right-hand side changes, but the matrix stays the same. An example of this scenario is solving time-dependent problems. It is possible then to start solving immediately, without deflation, then using the vectors resulting from the iterations to extract some spectral information which can be subsequently used to deflate the next linear systems, leading to a faster convergence. Below, we show how this second strategy of finding deflation vectors can be used in our case, by solving saddle point problems with CRAIG \cite{arioli2013generalized}.
	 
	 Our approach is inspired by the use of Ritz values and vectors associated with CG, more precisely with the Lanczos process with which CG is connected. The Ritz values are the eigenvalues of the tridiagonal matrix that can be constructed using the scalars computed by CG, and are approximations of the eigenvalues of the system matrix (see \cite[Chapter~5]{van2003iterative}). The  approximate eigenvectors associated with these approximate eigenvalues are known as Ritz vectors and they constitute the basis used in the deflation process. 
	 For deflation using Ritz vectors, harmonic Ritz vectors and others in the context of solving with GMRES, we refer to \cite{mercier2017new, sifuentes2013gmres, daas2021recycling}. 
	 
	 The particular strategy that we build on was presented in \cite{stathopoulos2010computing} for computing Ritz vectors from CG runs.
	 The key feature of their algorithm, called \textit{eigCG}, lies in the quick convergence of the Ritz vectors to the eigenvectors, without affecting the convergence of the linear solver. Even if an eigensolver and a linear solver may rely on the same algorithm (in this case, the Lanczos tridiagonalization), it is used differently in each case, which typically means one has to choose which objective gets prioritized: either the linear solver or the eigensolver. The authors of \cite{stathopoulos2010computing} show that their method can achieve both goals quickly. The success of their algorithm is partly due to the use of two consecutive tridiagonal matrices built from the output of CG. The mechanism allows them to obtain spectral information at a rate comparable to the more expensive full and reorthogonalized Lanczos process. The eigenvalue problem is solved relying on the linear solver, but not influencing it.

	 We follow the strategy given in \cite{stathopoulos2010computing} for computing Ritz vectors from CG runs and adapt it such that it fits our CRAIG solver and the underlying bidiagonalization process. The changes that we have to make stem from the fact that the generalized Golub-Kahan process works with two sets of orthogonal vectors (instead of one for CG) and the fact that we are looking for elliptic singular vectors, rather than eigenvectors. We give our interpretation in \Cref{alg:ESVD3}.

	 \begin{algorithm}[htbp] 
	 	\caption{Elliptic Singular Value Decomposition 3}
	 	\label{alg:ESVD3}
	 	\begin{algorithmic}[1] 	
	 		\REQUIRE \hspace{-0.2cm}:  $ \Am $ - rectangular matrix; $\quad$
	 		$ k $ - number of desired triplets;  $\quad$
	 		$ target $ - choose between $ smallest $ or $ largest $ triplets;  $\quad$	 		
	 		$ \Um_b, \Bm , \Vm_b $ from the bidiagonalization;  $\quad$
	 		$ \Um, \mSigma , \Vm$ previous approximate elliptic singular triplets.
	 		
	 		\ENSURE \hspace{-0.2cm}:   $ \Um, \mSigma , \Vm$ new approximate elliptic singular triplets. 
	 		\vspace{.5cm}

	 		\IF {$\Um, \mSigma, \Vm$ are given}
	 		\STATE $\vvv = $ first column of $\Vm _b; \quad \rv = \Um ^T \Am \vvv; \quad $
	 		$ \Bm =
	 		\begin{bsmallmatrix}
	 			\mSigma \ \mZ \\
	 			\mZ \ \Bm
	 		\end{bsmallmatrix} ;$		         $ \quad \Bm_{1:2k,2k+1}=\rv; $ 
	 		\STATE $ \Um _b =[ \Um \ \Um _b]; \quad \Vm _b =[ \Vm \ \Vm _b];  $
	 		\ENDIF
	 		
	 		\STATE Compute the SVD of $ \Bm  $ and keep the $ k $ $ target $ triplets in  $ \Um_{s1}, \mSigma_{s1}  $ and $ \Vm_{s1} $;
	 		
	 		\STATE Compute the SVD of $ \Bm_{1: \eta -1,1: \eta -1}  $ and keep the $ k $ $ target $ triplets in  $ \Um_{s2}, \mSigma_{s2}  $ and $ \Vm_{s2} $;
	 		
	 		\STATE Append a row of zeros to $ \Um_{s2}  $ and $ \Vm_{s2} $ so that they also have $\eta$ rows;
	 		
	 		\STATE Concatenate and orthogonalize the left vectors $ \Um_{so} = ortho ([\Um_{s1} \ \Um_{s2}]) $;
	 		\STATE Concatenate and orthogonalize the right vectors $ \Vm_{so} = ortho ([\Vm_{s1} \ \Vm_{s2}]) $;
	 		
	 		\STATE $ \Hm = \Um_{so} \add{^T} \Bm \Vm_{so} $;
	 		\STATE Compute the SVD of $ \Hm $ and store the triplets in  $ \Um_{s3}, \mSigma_{s3}  $ and $ \Vm_{s3} $;
	 		
	 		\STATE $ \Um = \Um_b \Um_{so} \Um_{s3}  $;
	 		$\quad$ $ \Vm = \Vm_b \Vm_{so} \Vm_{s3}  $;
	 		$\quad$ $ \mSigma = \mSigma_{s3}  $;

	 	\end{algorithmic}
	 \end{algorithm}

	 We describe the CRAIG solver in \Cref{sec:appendix} and \Cref{alg:CRAIG}. In particular, we highlight where \Cref{alg:ESVD3} can be integrated in order to obtain approximations to the elliptic singular triplets from the information computed by the solver. \Cref{alg:ESVD3} is to be called repetitively and proceeds in the following way. The first time, we have $\Bm \in \bR^{\eta \times \eta}$ with $\Um _b \in \bR^{m \times \eta}$ and $\Vm _b \in \bR^{n \times \eta}$ from the bidiagonalization, similar to the first step of \Cref{alg:ESVD2}. Three SVDs follow, with their results subscripted by $s1$, $s2$ and $s3$ respectively. Since the SVDs operate on small matrices, they are not particularly expensive. The first SVD is applied to $\Bm$, the second to $\Bm_{1: \eta -1,1: \eta -1} $, and the third to the combined result of the previous two. For this combination to be possible, $ \Um_{s2}  $ and $ \Vm_{s2} $ need to be extended and orthogonalized to $ \Um_{s1}  $ and $ \Vm_{s1} $. Collecting the results of the three SVDs, we have $ \Um, \mSigma , \Vm,$  a set of $2k$ approximate elliptic singular triplets, which concludes one step of \Cref{alg:ESVD3}. Subsequent steps refine these $2k$ triplets from an additional $\eta - 2k$ provided by the bidiagonalization. Once  the CRAIG solver concludes, it returns the $k$ triplets that we are interested in along with the solution of the saddle point system.

	 There are a few notable differences between \Cref{alg:ESVD2} and \Cref{alg:ESVD3}. The latter makes use of the vectors resulting from the GK bidiagonalization (part of the CRAIG solver), but does not influence it at all, since it does not include an augmentation. The bidiagonalization runs only once, without restarts, in contrast to the way \Cref{alg:ESVD2} uses it. 
	 Another point is that instead of keeping only $k$ approximations using the most recent bidiagonal matrix,  \Cref{alg:ESVD3} keeps an extra set of $k$ from the previous one. Combining these and refining them with further $\eta -2k$ vectors from the bidiagonalization leads to rapid convergence of the spectral information.
	 
	 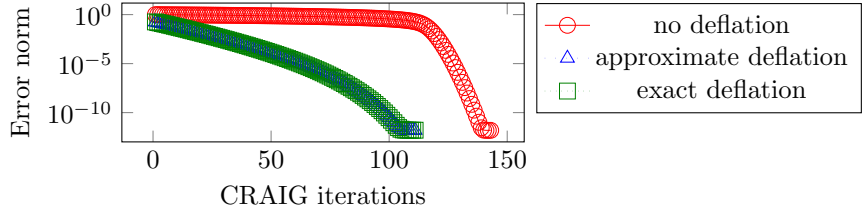
\begin{figure}[htbp]
    \centering
    \begin{tikzpicture}
    \begin{axis}[
    ymode=log,
    legend pos=outer north east,
	table/col sep=tab,
    width=\FigWidSolo \textwidth,
	height=\FigHeiSolo \textwidth,
    xlabel={CRAIG iterations},
     ylabel={Error norm}
    ]
	\addplot+[] table [x=no deflationX, y=no deflationY]{images/csv/conv1DCh512DeflApx5VsEx.csv}; 
\label{fig:item:conv1DCh512DeflApx5VsEx_no deflation} 
\addlegendentry{no deflation} 

\addplot+[] table [x=approximate deflationX, y=approximate deflationY]{images/csv/conv1DCh512DeflApx5VsEx.csv}; 
\label{fig:item:conv1DCh512DeflApx5VsEx_approximate deflation} 
\addlegendentry{approximate deflation} 

\addplot+[] table [x=exact deflationX, y=exact deflationY]{images/csv/conv1DCh512DeflApx5VsEx.csv}; 
\label{fig:item:conv1DCh512DeflApx5VsEx_exact deflation} 
\addlegendentry{exact deflation} 

\end{axis} 
\end{tikzpicture} 
\caption{Convergence curves of the CRAIG algorithm for the deflated 1D channel problem. Y-Axis: Energy norm of the relative error (Velocity field).} 
\label{fig:conv1DCh512DeflApx5VsEx} 
\end{figure}  
	 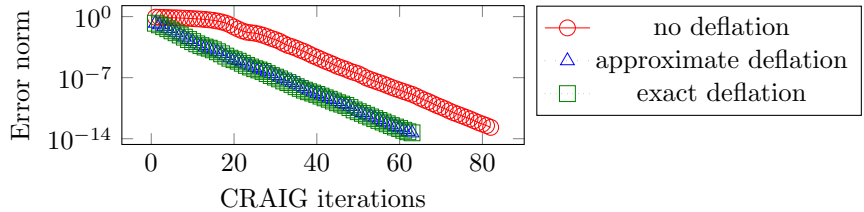
\begin{figure}[htbp]
    \centering
    \begin{tikzpicture}
    \begin{axis}[
    ymode=log,
    legend pos=outer north east,
	table/col sep=tab,
    width=\FigWidSolo \textwidth,
	height=\FigHeiSolo \textwidth,
    xlabel={CRAIG iterations},
     ylabel={Error norm}
    ]
	\addplot+[] table [x=no deflationX, y=no deflationY]{images/csv/convIFISS20DeflApxVsEx.csv}; 
\label{fig:item:convIFISS20DeflApxVsEx_no deflation} 
\addlegendentry{no deflation} 

\addplot+[] table [x=approximate deflationX, y=approximate deflationY]{images/csv/convIFISS20DeflApxVsEx.csv}; 
\label{fig:item:convIFISS20DeflApxVsEx_approximate deflation} 
\addlegendentry{approximate deflation} 

\addplot+[] table [x=exact deflationX, y=exact deflationY]{images/csv/convIFISS20DeflApxVsEx.csv}; 
\label{fig:item:convIFISS20DeflApxVsEx_exact deflation} 
\addlegendentry{exact deflation} 

\end{axis} 
\end{tikzpicture} 
\caption{Convergence curves of the CRAIG algorithm for the deflated IFISS channel problem.Y-Axis: Energy norm of the relative error (Velocity field).} 
\label{fig:convIFISS20DeflApxVsEx} 
\end{figure} 

	 In \Cref{fig:conv1DCh512DeflApx5VsEx} and \Cref{fig:convIFISS20DeflApxVsEx}, we plot the convergence history for the two test cases described in \Cref{sec:testPb} in the following three scenarios: no deflation, deflation using approximate vectors and deflation using exact vectors. The approximations were obtained using \Cref{alg:ESVD3} based on vectors coming from a previous linear solve of the same system. We can see how, in both cases, approximate vectors are just as effective as their exact counterparts in the process of deflation.

	 In \Cref{tab:ritzrr} we give the errors for the approximate elliptic singular values and vectors. If we denote the exact vectors by $\vvv^*_i$ and the approximate ones by $\Vm$, the quantities we report are $\normEu{\vvv_i^* - \Vm \Vm^T \vvv_i^*}$, with  $ |(\sigma _i-\sigma ^*_i)/ \sigma ^*_i|$ being the relative error for the values. We can see how recycling information from a linear solve can lead to fairly precise approximations. However, it may be the case that this level of accuracy cannot be reached after just one solve. This is especially true when attempting to identify many vectors. For this scenario, the authors of \cite{stathopoulos2010computing} present an algorithm  to  refine the set using information from several linear solves.
	 
	 \begin{table}
	 	\centering
	 	\caption{Relative errors for the smallest five approximate elliptic singular vectors and values (ordered descendingly) obtained from \Cref{alg:ESVD3}.}
	 	\label{tab:ritzrr}
	 	\begin{tabular}{cccccc}
	 		triplet index & 5th & 4th & 3rd & 2nd & 1st  \\
	 		\hline
	 		
	 		\multicolumn{6}{c}{1D Channel problem} \\
	 		vectors & 5.18E-06 & 3.52E-07 & 7.56E-08 & 3.23E-08 & 1.17E-08 \\
	 		values  & 3.78E-11 & 4.07E-13 & 4.19E-14 & 1.48E-15 & 1.29E-13  \\
	 		\hline
	 		
	 		\multicolumn{6}{c}{IFISS Channel problem}\\
	 		vectors & 4.56E-02 & 5.95E-05 & 7.02E-09 & 5.64E-08 & 1.02E-05            \\
	 		values  & 1.52E-04 & 5.37E-10 & 4.12E-15 & 5.94E-15 & 1.23E-13          
	 	\end{tabular}%
	 \end{table}

	 \section{Comparison with deflated MINRES}\label{sec:minres}
	 A popular choice for solving saddle point problems is the MINRES algorithm, applied monolithically and paired with an appropriately defined deflation to reduce convergence plateaus. In this section, we are interested in comparing the number of iterations performed by the linear solver following deflation, considering the two approaches: monolithic and segregated, with their respective solver MINRES and CRAIG.

	 In \cite{arioli2013generalized}, an important link between the generalized CRAIG algorithm and MINRES has been drawn. If the system in \Cref{eq:sad_po} is first preconditioned by a block-diagonal matrix with $\Wm$ and $\Id$ and then solved using MINRES, the required number of iterations will be double compared to that of generalized CRAIG. This is based on existing results regarding unpreconditioned MINRES and the CRAIG algorithm without the generalization, such as those in \cite{fischer2011polynomial} (Theorem 6.9.9). The same theorem highlights a situation where MINRES only reduces the residual every two iterations, by producing the same iterate twice in a row. This phenomenon occurs when the spectrum of the matrix is symmetric. The preconditioned saddle point matrix is 
	 \begin{equation*}
	 	\left[
	 	\begin{matrix}
	 		\Lm ^{-1} & \mZ \\
	 		\mZ & \Id
	 	\end{matrix}
	 	\right]
	 	\left[
	 	\begin{matrix}
	 		\Wm & \Am \\
	 		\Am ^{T} & \mZ 
	 	\end{matrix}
	 	\right]
	 	\left[
	 	\begin{matrix}
	 		\Lm ^{-T} & \mZ \\
	 		\mZ & \Id
	 	\end{matrix}
	 	\right]
	 	=
	 	\left[
	 	\begin{matrix}
	 		\Id& \tilde \Am \\
	 		\tilde \Am ^T  & \mZ
	 	\end{matrix}
	 	\right],
	 \end{equation*}
	 where we use a Cholesky decomposition of the (1,1)-block $ \Wm = \Lm \Lm ^T$, $ \Wm ^{-1}= \Lm ^{-T} \Lm ^{-1} $ and define $\tilde \Am = \Lm ^{-1} \Am $ as in \Cref{sec:defldefi}. We can describe its spectrum by adapting Result 1 from \cite{saunders1995solution} to our context:
	 \begin{equation*}
	 	\begin{cases}
	 		1, \quad &$m-n$ \text { times,} \\
	 		\frac{1}{2} + \sqrt{\sigma _i ^2 + \frac{1}{4}}, \quad &$n$ \text { times,} \\
	 		\frac{1}{2} - \sqrt{\sigma _i ^2 + \frac{1}{4}}, \quad &$n$ \text { times,} 
	 	\end{cases}
	 \end{equation*}
	 where $\sigma _i$ are the $n$ elliptic singular values of $\Am$ (singular values of $\tilde \Am $). If we ignore the cluster of values equal to one, which do not impact the convergence of MINRES,  the rest of the spectrum is symmetric around the value $0.5$. This symmetry then leads to redundant iterations in MINRES and GMRES as well, as proven theoretically in \cite{fischer1998minimum}. The authors do note that in numerical experiments, the two consecutive iterates may not necessarily be identical, but are very close nonetheless.   
	 
	 We can observe this effect in our context as well, when comparing generalized CRAIG and preconditioned MINRES. In \Cref{fig:1DCh_MINRESvsCRAIG}, we show the convergence curves of the two algorithms when solving the 1D channel problem described in \Cref{sec:testPb}. Notice how every second MINRES iteration does not visibly contribute to the reduction of the error. With the exception of this feature, the two algorithms converge in the same way, featuring the same two phases of slow initial convergence (plateau) and subsequent superlinear convergence.
	 
	 We proceed to test whether this 1:2 iteration ratio is preserved when applying deflation. When using MINRES, we consider the saddle point matrix as a whole and deflate it using the procedure described in \cite{gaul2013framework}, relying on the eigenpairs. When using CRAIG, we make use of the elliptic singular triplets and apply deflation following \Cref{sec:defldefi}. Both approaches eliminate the same number of values from the spectrum of their respective problem. The convergence history of both algorithms is given in {\Cref{fig:1DCh_MINRESvsCRAIG}}, for the 1D channel test case. We also compare the two scenarios on the IFISS channel problem presented in \Cref{sec:testPb}, with the results plotted in \Cref{fig:IFISS_MINRESvsCRAIG}. As we can see, the same relationship between the two algorithms holds, both with and without deflation. It should be noted that the respective curves are not identical (aside from the 1:2 scale) partly because the two algorithms under consideration have different minimization goals. For MINRES, this is the residual norm of the preconditioned saddle point system, whereas CRAIG minimizes the energy norm of the error for the primal variable $\uv$. 
	 
	 \begin{figure}
    \centering
	\begin{tikzpicture}
    \begin{groupplot}[
        group style={
            group size=2 by 1,
            horizontal sep=0pt,
            y descriptions at=edge left,
        },
		ymode=log,
		table/col sep=tab,
		width= \FigWidHalf \textwidth,
		height= \FigHeiHalf \textwidth,
    ]
    \nextgroupplot[xlabel={\ac{GKB} iterations}]
		\addplot+[] table [x=deflated 10 smallest valuesX, y=deflated 10 smallest valuesY]{images/csv/1DCh_CRAIG.csv}; 
		\label{fig:item:1DCh_CRAIG_deflated 10 smallest values} 
		
		\addplot+[] table [x=no deflationX, y=no deflationY]{images/csv/1DCh_CRAIG.csv}; 
		\label{fig:item:1DCh_CRAIG_no deflation} 
		
    \nextgroupplot[xlabel={MINRES iterations}]
        \addplot+[] table [x=deflated 10 smallest valuesX,     y=deflated 10 smallest valuesY]{images/csv/1DCh_MINRES.csv};
		\label{fig:item:1DCh_MINRES_deflated 10 smallest values}

		\addplot+[] table [x=no deflationX, y=no deflationY]{images/csv/1DCh_MINRES.csv}; 
		\label{fig:item:1DCh_MINRES_no deflation} 
    \end{groupplot}
    \node (title) at ($(group c1r1.center)!0.5!(group c2r1.center)+(0,2cm)$) {Energy norm of relative error (Velocity field)};
\end{tikzpicture}
\caption{Convergence curves of the CRAIG algorithm (Left) and preconditioned MINRES (Right) for the 1D channel problem. Undeflated in \ref{fig:item:1DCh_MINRES_no deflation}, exact deflation in \ref{fig:item:1DCh_MINRES_deflated 10 smallest values}. } 
\label{fig:1DCh_MINRESvsCRAIG} 
\end{figure}
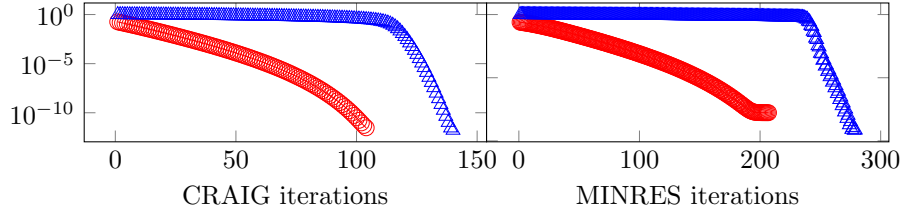 
	 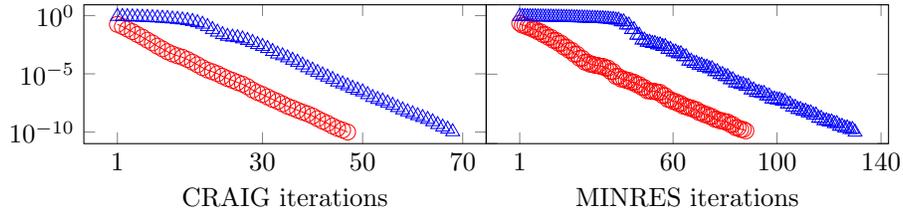
\begin{figure}
    \centering
	\begin{tikzpicture}
    \begin{groupplot}[
        group style={
            group size=2 by 1,
            horizontal sep=0pt,
            y descriptions at=edge left,
        },
		ymode=log,
		table/col sep=tab,
		width= \FigWidHalf \textwidth,
		height= \FigHeiHalf \textwidth,
    ]
    \nextgroupplot[xlabel={\ac{GKB} iterations},xtick={1,30,50,70}]
		\addplot+[] table [x=deflated 5 smallest valuesX, y=deflated 5 smallest valuesY]{images/csv/IFISS_CRAIG.csv}; 
		\label{fig:item:IFISS_CRAIG_deflated 5 smallest values}  
		
		\addplot+[] table [x=no deflationX, y=no deflationY]{images/csv/IFISS_CRAIG.csv}; 
		\label{fig:item:IFISS_CRAIG_no deflation} 
		
    \nextgroupplot[xlabel={MINRES iterations},xtick={1,60,100,140}]
        \addplot+[] table [x=deflated 5 smallest valuesX, y=deflated 5 smallest valuesY]{images/csv/IFISS_MINRES.csv}; 
		\label{fig:item:IFISS_MINRES_deflated 5 smallest values} 

		\addplot+[] table [x=no deflationX, y=no deflationY]{images/csv/IFISS_MINRES.csv}; 
		\label{fig:item:IFISS_MINRES_no deflation} 
    \end{groupplot}
    \node (title) at ($(group c1r1.center)!0.5!(group c2r1.center)+(0,2cm)$) {Energy norm of relative error (Velocity field)};
\end{tikzpicture}
\caption{Convergence curves of the CRAIG algorithm (Left) and preconditioned MINRES (Right) for the IFISS channel problem. Undeflated in \ref{fig:item:IFISS_MINRES_no deflation}, exact deflation in \ref{fig:item:IFISS_MINRES_deflated 5 smallest values}. } 
\label{fig:IFISS_MINRESvsCRAIG} 
\end{figure} 

	 \subsection{Augmented CRAIG}
	 Above, we have compared the two algorithms only with respect to the number of iterations. If we are interested in memory use or time, we can consider augmented MINRES \cite{olshanskii2010acquired} and augmented CRAIG, the latter being described below.

	 Augmentation can be seen as an efficient implementation of deflation. 
	 First, given the definition of the operators $\Mm$, $\Pm$ and $\Qm$ in \Cref{eq:deflM,eq:deflP,eq:deflQ}, they will most likely be dense matrices if explicitly computed and stored. 
	 Then, the deflated matrix, whether computed  as $\Pm\Am$ or $\Am\Qm$, will also be dense. However, in many applications, the matrices that result from the discretization are (very) large and sparse. Losing sparsity by a naive deflation might make the system difficult or impossible to solve due to excessive requirements concerning memory and time. 
	 
	 We have defined deflation to be paired with the CRAIG solver \cite{arioli2013generalized}, which only uses the matrix $\Am$ as part of matrix-vector products $\Am \xv$, for some vector $\xv$. To avoid the difficulties mentioned above, these products should be reformulated to include the effect of deflation in an implicit manner. Applying the definition of the deflation operators we get $\Am \xv - \Am \Vm \mSigma ^{-1} \Um ^T \Am \xv$, which has the advantage of manipulating matrices that are either sparse or small (with only $k$ columns). Suitable alternatives in the same style can be defined to replace the operations in \Cref{eq:simpler_corr}, avoiding explicit use of $ \Pm, \Qm, $ and $\Mm$. 
	 This strategy corresponds to an \textit{augmented} CRAIG solver, similar to the augmented MINRES described in  \cite{olshanskii2010acquired}, where the basis for deflation is passed to the solver and applied whenever needed. Here, we note that MINRES operates monolithically, with a basis $\Ym \in \bR ^{(m+n) \times k}$ applied repeatedly, whereas augmented CRAIG uses either $\Um \in \bR ^{m \times k}$ or $\Vm \in \bR ^{n \times k}$. Thus, adding augmentation to CRAIG is less expensive than to MINRES. Although CRAIG needs matrix-vector products with both $\Am$ and $\Am  ^T$, it is sufficient to reformulate only that with $\Am$. It can be shown using the properties from \Cref{eq:proj_prop} that the Schur complement $ \Sm = \Am ^T \Wm ^{-1} \Am $ is correctly deflated even with this "one-sided" augmentation, which leads CRAIG to converge quickly to the solution of \Cref{eq:simpler_defl}. For a more detailed comparison of augmentation and deflation, we refer to \cite{coulaud2013deflation,gaul2011deflated}.
	 
	 \begin{remark}
	 	Further savings compared to augmented MINRES are possible if we have exact elliptic singular triplets. Then, it is true that $\Vm \Vm ^T = \Mm \Am$, which we use to define a simpler and quicker matrix-vector product: $\Am \xv - \Am \Vm \Vm ^T \xv$, circumventing matrices of order $m$. This simplification can still be applied if we have approximate instead of exact triplets, provided we are only interested in a solution of moderate accuracy, as is the case in certain applications. One such example is solving linear systems as part of a nonlinear solver, as shown in \cite{daas2021recycling}, where the linear solver needs to reduce the residual norm only by $10^{-1}$. Using the simpler $\Qm := \Id_n - \Vm \Vm^T$ in the deflation and correction stages we described in \Cref{sec:defldefi} introduces a perturbation proportional to the accuracy of the elliptic singular vectors. This then limits the accuracy of the solution we seek for the initial saddle point system. Such a behavior is similar to a different setting, where we approximate the action of $\Wm ^{-1}$ in the CRAIG solver by using an iterative solver instead \cite{darrigrand2022inexact}.  
	 \end{remark}

	 \section{Deflation and preconditioning}
	 \label{sec:deflPrec}
	 
	 We have seen how targeted, local changes in the spectrum can be induced by means of deflation. Preconditioning, on the other hand, typically aims to change the spectrum as a whole, clustering the values away from zero. For saddle point problems, the following block diagonal preconditioner is considered ideal 
	 	\begin{equation*}
	 		\cP = 
	 		\left[
	 		\begin{matrix}
	 			\Wm & \mZ \\
	 			\mZ & \Sm
	 		\end{matrix}
	 		\right].
	 	\end{equation*}
	 	It is ideal in the sense that depending on the solver used, when applied as $\cP ^{-1}$ to a system like the one in \Cref{eq:sad_po}, three or fewer iterations  will be sufficient to converge to the exact solution. 
	 	Since $\Sm ^{-1}$ is difficult to obtain, an approximation is used instead, trying to balance the cost of building and applying it against the improvement in convergence speed that it brings. The preconditioning strategy used in the remainder of this section, but also other choices for saddle point problems are described in \cite[Section~10]{benzi2005numerical}.
	 
	 For Stokes flow problems, the pressure mass matrix $\Qm _p$ is commonly employed to achieve such an approximation, exploiting the spectral equivalence it has to the Schur complement (see \cite[Chapter~3]{elman2014finite}). It is formally defined as 
	 	\begin{equation*}
	 		\Qm _p=[\qv _{ij}], \quad \qv _{ij}= \int _{\Omega}  \psi _i \psi _j, \quad i,j=1,...,n,
	 	\end{equation*}
	 	where $\Omega$ is the domain, $\{ \psi \} _{i=1} ^{n}$ are basis functions for the finite-dimensional space $\cM _n \in \cM$, with $\cM$ being the Hilbert space for the pressure. 
	 To further reduce the cost of applying the preconditioner, it is simplified to a diagonal matrix, for example by lumping, i.e. placing each row-sum on the diagonal. 
	 
	 We proceed to apply this kind of preconditioning to our IFISS channel problem and compare its behavior with that of deflation. As \Cref{fig:IFISS_PrecDefl} illustrates, CRAIG converges faster in the presence of this preconditioner, but the initial plateau phase is mostly unaffected. Deflation of the five smallest elliptic singular values removes the plateau, enabling a faster convergence to solutions of moderate accuracy (error norm $>10^{-7}$). However, this deflation does not improve the second stage of the convergence process, which is important if we aim for a very accurate solution. The complementary action of the preconditioner and deflation motivates the third scenario, in which both are applied, leading to the best result. A similar combined approach was employed in \cite{olshanskii2010acquired} for a different type of Stokes problem solved with MINRES.

	 \begin{figure}[h]	 	
	 		\centering	
	 		\begin{subfigure}{\FigWidHalf \textwidth}
    \centering
    \begin{tikzpicture}
    \begin{axis}[
	table/col sep=tab,
    ymode=log,    
    legend pos=north east,
    width=\FigWid \textwidth,
	height=\FigHei \textwidth,
    xlabel={CRAIG iterations},
    ylabel={Error norm}
    ]
	
\addplot+[] table [x=originalX, y=originalY]{images/csv/IFISS_PrecDefl.csv}; 
\label{fig:item:IFISS_PrecDefl_original} 

\addplot+[] table [x=deflatedX, y=deflatedY]{images/csv/IFISS_PrecDefl.csv}; 
\label{fig:item:IFISS_PrecDefl_deflated} 

\addplot+[] table [x=preconditionedX, y=preconditionedY]{images/csv/IFISS_PrecDefl.csv}; 
\label{fig:item:IFISS_PrecDefl_preconditioned} 

\addplot+[] table [x=defl. & prec.X, y=defl. & prec.Y]{images/csv/IFISS_PrecDefl.csv}; 
\label{fig:item:IFISS_PrecDefl_defl. & prec.} 

\end{axis} 
\end{tikzpicture} 
\caption{CRAIG convergence curves. Y-Axis: Energy norm of the relative error (Velocity field)} 
\label{fig:IFISS_PrecDefl} 
	 		\end{subfigure}
	 		\begin{subfigure}{\FigWidHalf \textwidth}
    \centering
    \begin{tikzpicture}
    \begin{axis}[
	table/col sep=tab,
    ymode=log,    
    legend pos=north east,
    width=\FigWid \textwidth,
	height=\FigHei \textwidth,
    xlabel={Ell. sing. value. index},
    ylabel={Magnitude}
    ]
	
\addplot+[] table [x=originalX, y=originalY]{images/csv/IFISS_spec_PrecDefl.csv}; 
\label{fig:item:IFISS_spec_PrecDefl_original}

\addplot+[] table [x=deflatedX, y=deflatedY]{images/csv/IFISS_spec_PrecDefl.csv}; 
\label{fig:item:IFISS_spec_PrecDefl_deflated}

\addplot+[] table [x=preconditionedX, y=preconditionedY]{images/csv/IFISS_spec_PrecDefl.csv}; 
\label{fig:item:IFISS_spec_PrecDefl_preconditioned} 

\addplot+[] table [x=defl. & prec.X, y=defl. & prec.Y]{images/csv/IFISS_spec_PrecDefl.csv}; 
\label{fig:item:IFISS_spec_PrecDefl_defl. & prec.} 

\end{axis} 
\end{tikzpicture} 
\caption{Distribution of elliptic singular values.} 
\label{fig:IFISS_spec_PrecDefl} 
	 		\end{subfigure}

	 	\caption{ Corresponding CRAIG convergence curves and Schur complement spectra for the IFISS channel of length 20. We consider the original case \ref{fig:item:IFISS_PrecDefl_original}, the deflated case (five smallest elliptic singular values) \ref{fig:item:IFISS_PrecDefl_deflated}, the preconditioned case (lumped pressure mass matrix) \ref{fig:item:IFISS_PrecDefl_preconditioned}, and the case with both deflation and preconditioning \ref{fig:item:IFISS_PrecDefl_defl. & prec.}. }
	 	
	 \end{figure}

	 To understand the effect of the pressure mass matrix preconditioner, we plot in \Cref{fig:IFISS_spec_PrecDefl} the elliptic spectra for the four cases illustrated in \Cref{fig:IFISS_PrecDefl}. The most visible effect is how almost 600 of the largest values cluster around 1, while a second cluster forms around 0.7. However, the small outliers are not brought closer to the rest of the values, which is reflected in the persistence of the plateau visible in \Cref{fig:IFISS_PrecDefl}.
	 
	 Stokes problems on elongated domains represent a special case, where the spectrum of the Schur complement changes also as a function of the anisotropy. More precisely, the smallest elliptic singular value is bounded (see \cite{chizhonkov2000domain}) 
	 	\begin{equation*}
	 		\frac{1}{2 \sqrt{15}} l^{-1} \leq \sigma _n \leq \frac{\pi}{2 \sqrt{3}} l^{-1},
	 	\end{equation*}
	 	where $l>1$ is the ratio between the sides of the channel. As this ratio increases, the smallest elliptic singular value approaches zero. Additionally, as seen in \Cref{fig:EllSpecIFISSChL}, more and more small outliers appear. Deflation can then be used to manage the outliers' effect, while the pressure mass matrix as a preconditioner clusters the other elliptic singular values.

     \section{Conclusions}
     \label{sec:conc}
     
     In this paper, we have defined a particular kind of deflation to be used with CRAIG, a segregated linear solver for symmetric saddle point problems. The motivation for this development was to accelerate convergence for problems with a spectrum which contains outliers, leading to stagnation periods in the convergence process. As an example of such a problem, we have considered two kinds of Stokes flow in thin, elongated channels. 
     Our definition of deflation uses bases from an elliptic singular value decomposition of the off-diagonal block $\Am$ as an alternative to explicitly deflating the Schur complement, which could pose numerical difficulties. With this deflation mechanism, stagnation periods can be shortened or removed, even if we only have bases with inexact spectral information.
     We have extended two iterative algorithms, such that their output is a partial approximate elliptic singular value decomposition which can be employed in deflation. 
     Following a comparison with MINRES, we have seen that our deflation for a segregated solver performs just as well as that for a monolithic one, removing the plateau in convergence. The ability to exploit block structure and the lack of redundant iterations in CRAIG highlights the advantage and need of having a segregated deflation strategy. Additionally, the augmentation procedure we have described for CRAIG operates with smaller matrices than MINRES, leading to faster iterations.
     As future work, a promising avenue is that of using multigrid transfer operators instead of the ESVD as a way to obtain bases for deflation.

	 \bibliographystyle{siamplain}
	 \bibliography{references}
     
     \pagebreak
     \appendix
     \section{The CRAIG algorithm}
     \label{sec:appendix}
     
     The CRAIG algorithm iteratively solves the saddle point system in \Cref{eq:sad_po} in the following way. It computes the generalized Golub-Kahan bidiagonalization
     \begin{equation}
     	\label{eq:oriGKB}    
     	\Am \Vm _i = \Wm \Um _i \Bm _i, \quad  \Um _i ^T \Wm \Um _i = \Id _i, \quad
     	\Am ^T \Um _i = \Vm _i \Bm ^T _i + \beta _{i+1} \qv _{i+1} \ev _i^T,  \quad \Vm _i ^T  \Vm _i = \Id _i,    
     \end{equation}
     with the bidiagonal matrix
     \begin{equation}
     	\label{eq:BmatGKB}
     	\Bm _i=
     	\left[ 
     	\begin{matrix}
     		\alpha_1 & \beta_2 & 0 & \ldots & 0 \\
     		0 &	\alpha_2 & \beta_3 &  \ldots & 0 \\
     		\vdots &	\vdots & \vdots & \vdots & \vdots  \\
     		0 & \ldots & 0 & \alpha_{i-1} & \beta_i \\
     		0 & \ldots & 0 & 0 & \alpha_{i}
     	\end{matrix}
     	\right].
     \end{equation}
     Using $\Um _i$ and $\Vm _i$ it transforms \Cref{eq:sad_po} into
     \begin{equation}\label{eq:transfSys}
     	\left[
     	\begin{matrix}
     		\Id  _i& \Bm _i   \\
     		\Bm _i  ^T & \mZ 
     	\end{matrix}
     	\right]
     	\left[
     	\begin{matrix}
     		\zv _i   \\
     		\yv  _i
     	\end{matrix}
     	\right]
     	=
     	\left[
     	\begin{matrix}
     		\mZ  \\
     		\Vm _i ^T \bv 
     	\end{matrix}
     	\right], 
     \end{equation}
     with the solution
     \begin{equation}
     	\label{eq:zandy}
     	\zv _i= \beta _1  \Bm _i ^{-T} \ev _1; \quad \yv _i= - \Bm _i ^{-1} \zv _i,
     \end{equation}
     where $\ev _1 $ is the first column of the identity matrix $\Id _i \in \bR ^{i \times i}$. The approximate solution of the initial problem is then given by
     \begin{equation}
     	\label{eq:GKBapx}
     	\uv ^{(i)} = \Um _i \zv _i; \quad \pv ^{(i)} = \Vm _i \yv _i.
     \end{equation}
     Note: $\uv ^{(i)} $ is the $i$-th approximate solution, while $\uv _i $ is the $i$-th left vector.

     It is also possible to define updates for $\uv ^{(i)} $ and $\pv ^{(i)} $ in a recursive manner. To this end, we will need only $\zeta _i$ (the last entry of $\zv _i$) and an auxiliary vector $\dv$ for $\pv ^{(i)} $ (see \cite{arioli2013generalized} for details). The advantage of a recursive approach is that each update  needs only the latest pair of left and right vectors.

     In \Cref{alg:CRAIG}, we give a description of the CRAIG solver, and how it can interact with \Cref{alg:ESVD3}. An efficient implementation has been given in \cite{arioli2013generalized}, making use of the simple structure of $\Bm$ and short term recurrences. The author also presents error bounds useful for monitoring convergence.  
     
     An important point is that the use of short term recurrences implies there is no need to store the left and right vectors in order to compute the approximate solutions. However, for our purposes in this work they are necessary in order to extract information about the elliptic singular triplets from them using \Cref{alg:ESVD3}. After a group of $\eta$ left and right bidiagonalization vectors have been processed, they can be discarded.
     Additionally, while the iterative linear solver can deliver satisfactory approximate solutions even if the bidiagonalization vectors gradually lose orthogonality, the same may not necessarily hold in our extended context. This leads to the potential need for reorthogonalization, with similar remarks as in \Cref{sec:compESVD}.

     We present in \Cref{alg:CRAIG} how the original Craig solver is extended to make use of \Cref{alg:ESVD3}. We mark the additional statements needed for \Cref{alg:ESVD3} with special line numbers, and the associated inputs and outputs with a  $\bullet$.
     
     \begin{algorithm}
     	\caption{CRAIG algorithm}
     	\label{alg:CRAIG}
     		\begin{algorithmic}[1]
     			\REQUIRE \hspace{-0.2cm}:  
     			$ \Am \in \bR^{m \times n}, (n<m) $; $\quad$
     			$ \Wm \in \bR^{m\times m}, $ SPD; $\quad$
     			$ \gvv \in \bR^{m} $; $\quad$
     			$ \rv  \in \bR^{n} $; $\quad$
     			$ maxiter  \in \bR $; $\quad$ 
     			\\ $ \bf{ \bullet \quad} $ $ k  \in \bR $ - number of elliptic singular triplets; $\quad$
     			\\ $ \bf{ \bullet \quad} $ $ \eta  \in \bR $ - maximum search subspace dimension; $\quad$
     			\\ $ \bf{ \bullet \quad} $ $ \textit{target} $ - choose between smallest or largest triplets;
     			\ENSURE \hspace{-0.2cm}:
     			$ \uv \in \bR^{m} $; $\quad$
     			$ \pv  \in \bR^{n} $; $\quad$
     			\\ $ \bf{ \bullet \quad} $ $ \Um  \in \bR^{m \times k} $ - approximate left elliptic singular vectors; $\quad$
     			\\ $ \bf{ \bullet \quad} $ $ \mSigma  \in \bR^{k \times k} $ - approximate elliptic singular values; $\quad$
     			\\ $ \bf{ \bullet \quad} $ $ \Vm  \in \bR^{n \times k} $ - approximate right elliptic singular vectors; $\quad$
     			\vspace{.5cm}  
     			\STATE{ $ \uv^{(0)} = \Wm ^{-1} \gvv; \quad$ $\quad  \bv = \rv - \Am ^T \uv^{(0)};  \quad
     				\beta_1 = \normEu{\bv}; \quad  \vvv_1 = \bv / \beta_1; $}
     			\STATE{$\wv = \Wm^{-1} \Am \vvv_1; \quad \alpha_1 = \normM{\wv}; \quad \uv_1 = \wv / \alpha_1; \quad  
     				\zeta _1 = \beta _1 / \alpha_1; \quad \dv _1 = \vvv _1 / \alpha_1; \quad $}    
     			\STATE{$ \uv^{(1)} = \uv^{(0)} + \zeta _1 \uv_1; \quad \pv^{(1)} = - \zeta _1 \dv_1;   $ $ \quad  i = 2;$}
     			\separateState{\STATE  $ \Bm _1 = \alpha _1; \quad \Vm _1 = \vvv _1; \quad \Um _1 = \uv _1;$}
     			\WHILE{not converged AND $i<maxiter$}
     			\STATE{$\vvv = \Am^T \uv_{i-1} - \alpha_{i-1} \vvv_{i-1} ;  $} 
     			\separateState{\STATE  Reorthogonalization: $\vvv= \vvv - \Vm _i \Vm _i ^T \vvv; $}
     			\STATE {$ \beta_{i} = \normEu{\vvv}; \quad \vvv_{i} = \vvv / {\beta_{i}}; $}
     			\STATE{   $ \wv = \Wm^{-1} \left(  \Am \vvv_{i} - \beta_{i} \Wm \uv_{i-1} \right);  $}
     			\separateState{\STATE{  $ \Wm $-reorthogonalization:	$ \wv= \wv - \Um _i \Um _i ^T \Wm \wv; $}}
     			\STATE{$ \alpha_{i} = \normM{\wv}; \quad  \uv_{i} = \wv / \alpha_{i}; $
     				
     				\separateState{\STATE 
     					$\Vm _{i}= [\Vm _{i-1} \ \vvv _{i}];  \quad  $ $ \Um _{i}  =[ \Um _{i-1}  \ \uv _{i}]; \quad \Bm _{i} = [\Bm _{i-1} \ \beta _{i} \ev_{i-1} ];  $
     					$ \quad$
     					$ \Bm _{i}= 
     					\begin{bsmallmatrix}
     						\Bm _{i} \\
     						\alpha _{i} \ev_{i}^T
     					\end{bsmallmatrix} ;$} }
     			
     			\updatelinenoprint
     			\addtocounter{ALC@line}{1} 
     			\IF{$mod(i,\eta) = 0$}
     			\STATE  $\Bm _{\eta} = $ lower right $\eta \times \eta$ submatrix of $\Bm _{i}$;
     			\IF{$i = \eta$}
     			\STATE  Call \Cref{alg:ESVD3} with $\Am, k, \Um _b = \Um _i, \Bm = \Bm _{\eta} ,\Vm _b = \Vm _i. $ Get $2k$ approximate triplets in $\Um, \mSigma, \Vm$.
     			\STATE  $\eta = \eta -2 k;$
     			\ELSE
     			\STATE  Call \Cref{alg:ESVD3} with $\Am, k, \Um _b = \Um _i, \Bm = \Bm _{\eta} ,\Vm _b = \Vm _i$ and the current $ \Um, \mSigma, \Vm. $ Get $2k$ updated approximate triplets in $\Um, \mSigma, \Vm$.
     			\ENDIF
     			\STATE  Clear $\Vm _i$ and $\Um _i$.
     			\ENDIF
     			\revertlinenoprint
     			\STATE{ $  \zeta _{i} = - \frac{\beta _{i} }{\alpha _{i}} \zeta _{i-1} ; \quad \dv _{i} = (\vvv _{i}  - \beta_{i} \dv _{i-1}) / \alpha _{i};$    }
     			\STATE{$ \uv^{(i)} = \uv^{(i-1)} + \zeta _{i} \uv _{i}; \quad \pv^{(i)} = \pv^{(i-1)} - \zeta _{i} \dv _{i};  $}
     			\STATE{$i = i+1;$}
     			\ENDWHILE
     			\STATE{ $\uv = \uv^{(i)}; \quad \pv= \pv^{(i)}; \quad$ }
     			\separateState{\STATE{ Choose the $ k $ $ target $ triplets from  $ \Um, \mSigma  $ and $ \Vm $. }}
     		\end{algorithmic}     	
     \end{algorithm}

\end{document}